\documentclass[11pt,twoside]{article}
\usepackage{amssymb}

\usepackage{latexsym}
\usepackage{graphicx}
\usepackage{amsmath}
\usepackage{hyperref,amsthm}

\numberwithin{equation}{section}

\normalsize \makeatletter

\newcommand{\be}{\begin{equation}}
\newcommand{\ee}{\end{equation}}
\newcommand{\bea}{\begin{eqnarray}}
\newcommand{\eea}{\end{eqnarray}}

\newtheorem{thm}{Theorem}[section]
\newtheorem{cor}[thm]{Corollary}
\newtheorem{lem}[thm]{Lemma}
\newtheorem{prop}[thm]{Proposition}

\newtheorem{rem}[thm]{Remark}
\newtheorem{cl}[thm]{Claim}

\begin{document}

\title{$L^2$-concentration phenomenon for Zakharov system\\ below energy norm\thanks{This work is
   supported  by NSFC 10571158}}

\author{ Daoyuan Fang, Sijia Zhong\\ Department of Mathematics,
  Zhejiang University, \\
Hangzhou 310027, China }

\maketitle \baselineskip 17pt
\begin{abstract}
In this paper, we'll prove a $L^2$-concentration result of Zakharov
system in space dimension two, with radial initial data
$(u_0,n_0,n_1)\in H^s\times L^2\times H^{-1}$ ($\frac
{16}{17}<s<1$), when blow up of the solution happens by
$I$-method. In additional to that we find a blow up character of
this system. Furthermore, we improve the global existence result
of Bourgain's to above spaces.
\\\\
\textit{Keywords}: Zakharov system in space dimension two;
$L^2$-concentration; blow up; global existence
\end{abstract}


\section{Introduction}

In this paper, we consider the following Zakharov system in space
dimension two: \be\left\{
\begin{array}{ll}
iu_t+\Delta u=nu, \label{1.1}\\
\Box n=\partial_{tt}n-\triangle n=\triangle |u|^2, \\
u(0,x)=u_0(x),\ n(0,x)=n_0(x),\ n_t(0,x)=n_1(x),
\end{array}
\right.\ee where $\triangle$ is the Laplacian operator in
$\mathbb{R}^2$, $u:[0,T)\times\mathbb{R}^2\rightarrow\mathbb{C}$,
$n:[0,T)\times\mathbb{R}^2\rightarrow\mathbb{R}$, and $u_0$,
$n_0$, $n_1$ are the initial data. We consider the Hamiltonian
case, that is, we assume that there is a
$w_0:\mathbb{R}^2\rightarrow\mathbb{R}$ such that
$n_t(0)=n_1=-\triangle w_0$. Then for any $t$, there is a $w(t)$
such that $n_t(t)=-\triangle w(t)=-\nabla\cdot v(t)$, where
$v(t)=\nabla w(t)$. In this case, (\ref{1.1}) can be written in
the form \be\left \{\begin{array}{ll}
iu_t+\Delta u=nu, \label{1.2}\\
n_t=-\nabla\cdot v, \\
v_t=-\nabla n-\nabla|u|^2, \\
u(0,x)=u_0(x),\ n(0,x)=n_0(x),\ v(0,x)=v_0(x),
\end{array}
\right.\ee

The Zakharov system was introduced in \cite{Z} to describe the long
wave Langmuir turbulence in a plasma. The function $u$ represents
the slowly varying envelope of the rapidly oscillating electric
field, and the function $n$ denotes the deviation of the ion density
from its mean value. We usually place the initial data $u_0\in H^k$,
the initial position $n_0\in H^l$ and the initial velocity $n_1\in
H^{l-1}$ for some real $k$, $l$.

It is well-known that the Schr\"{o}dinger equation is invariant
under the dilation transformation \[u(t,x)\rightarrow
u_\lambda(t,x)=\lambda u(\lambda^2t,\lambda x),\] while the wave
equation is invariant with the following transformation
\[n(t,x)\rightarrow n_\lambda(t,x)=\lambda n(\lambda t,\lambda x).\]
However, the Zakharov system doesn't have a true scale invariance
because the two relevant dilation transformations are incompatible.
Nevertheless the critical regularity is $(k,l)=(-\frac{1}{2},0)$.

For the local existence theory about this system. From \cite{GTV},
one can see that when $d=2$,
 the Cauchy problem
(\ref{1.1}) with $(u_0,n_0,n_1)\in H^k\times H^l\times H^{l-1}$ is
local well posed if $l\geqslant 0$ and $2k-(l+1)\geqslant0$.
Therefore the lowest allowed values of $(k,l)$ is $(\frac{1}{2},0)$.

On the other hand, if we replace $\Box$ in (\ref{1.1}) by
$\Box_c=c^{-2}\partial_{tt}-\triangle$, i.e. introducing explicitly
the ion sound velocity, then considering the limit
$c\rightarrow\infty$, the system (\ref{1.1}) reduces formally to the
nonlinear Schr\"{o}dinger equation \be\label{1.3}iu_t+\triangle
u=-|u|^2u,\ee which is just the $L^2$-critical focusing case for
$d=2$.

As for this Schr\"{o}dinger equation, the results in \cite{H},
\cite{H2} and so on for $s=1$, and \cite{CRSW}, \cite{FZ}, for
$1>s>\frac{1+\sqrt{11}}{5}$, tell us that there is some
$L^2$-concentration phenomenon for finite time blow up solutions,
i.e.
\[\limsup_{t\uparrow T^*}
\sup_{\scriptsize{\begin{array}{c}B\subset \mathbb{R}^2\\
R(B)\leqslant(T^*-t)^{\frac{s}{2}-}\end{array}}}\int_B|u(t)|^2\geqslant\|Q\|_{L^2}^2.\]
Here, $Q$ is the ground state for Schr\"{o}dinger equation, that is,
the unique positive solution (up to translations) of \be\label{1.4}
\triangle Q-Q+|Q|^2Q=0.\ee

In \cite{AA2}, \cite{OT1} and \cite{SW} the convergence of the
solutions of the $c$ dependent Zakharov system to those of NLS
equation when $c\rightarrow\infty$ was studied, which implies that
the $L^2$-concentration phenomenon like $L^2$-critical focusing
Schr\"{o}dinger equations may also happen. Glangetas and Merle in
\cite{GM1} proved this phenomenon for $(k,l)=(1,0)$ which is the
energy case.

We are interested here in the $L^2$-concentration phenomenon for
$s<1$ when blow up occurs of Zakharov system as well. What we want
to show is for some $0<k<1$ this phenomenon also holds true:

\begin{thm}\label{th1}
For $(u_0,n_0,n_1)\in H^s\times L^2\times H^{-1}$, radial,
$\frac{16}{17}<s<1$, if $(u,n)$ is a blow-up solution to equation
(\ref{1.1}), i.e. $T^*<\infty$ is its maximum existing time, then
there is a constant $m_n>0$ depending on the initial data such that
the following properties are true: $\forall R>0$, \be\label{1.5}
\limsup_{t\rightarrow T^*} \|u(t,x)\|_{L^2(|x|\leqslant
R)}\geqslant\|Q\|_{L^2},\ee and \be\label{1.6} \limsup_{t\rightarrow
T^*}\|n(t,x)\|_{L^1(|x|\leqslant R)}\geqslant m_n.\ee
\end{thm}
\begin{rem}
We can't remove the radial requirement because of the endpoint
Strichartz estimate for Schr\"{o}dinger equation we needed.
\end{rem}

As a quick result of the above theorem, and by the conservation of
$L^2$-norm of $u$, one has:

\begin{cor}\label{cor3.1}
For $(u_0,n_0,n_1)\in H^s\times L^2\times H^{-1}$, radial,
$\frac{16}{17}<s<1$, if $\|u_0\|_{L^2}<\|Q\|_{L^2}$, then the
corresponding solution to (\ref{1.1}) is global, i.e. $T^*=\infty$.
\end{cor}

In fact, the global well posedness for $k=l+1\geqslant3$ and small
data is considered in \cite{AA1}. Then Bourgain \cite{B1}
\cite{B2} introduced a new method to study the Cauchy problem for
nonlinear dispersive evolution equation, and applied it in
\cite{BC} to prove well posedness (both local and global) for
finite energy solutions namely for $k=l+1=1$ (also with small
initial data). Therefore, the above result is a  improvement of
the former result.

Now, let's briefly state about the proofs to Theorem \ref{th1}.

As we consider the Hamiltonian case, there are two conservations:
mass and energy (if exists). $\forall t\in[0,T^*)$, \be\label{1.7}
\int_{\mathbb{R}^2}|u(t,x)|^2dx=\int_{\mathbb{R}^2}|u_0(x)|^2dx,\ee
\be\label{1.8} H(t)=H(u(t),n(t),v(t))=H(u_0,n_0,v_0)=H_0,\ee where
\be\label{1.9}H(u,n,v)=\int_{\mathbb{R}^2}|\nabla
u(t,x)|^2+n(t,x)|u(t,x)|^2+\frac{1}{2}n^2(t,x)+\frac{1}{2}|v(t,x)|^2dx\ee
and $v$ has been defined before.

First, we split $n$ into its positive and negative frequency parts
according to \be\label{1.10} n_{\pm}=n\pm
i\Lambda^{-1}\partial_tn,\ee where $\Lambda=\sqrt{-\triangle}$. Thus
$n=\frac{n_++n_-}{2}$, $n_+=\bar{n}_-$, and equation (\ref{1.1})
equals to \be\label{1.11} \left\{\begin{array}{ll}
iu_t=-\triangle u+\frac{n_++n_-}{2}u \\
(i\partial_t\mp\Lambda)n_\pm=\mp\Lambda^{-1}\Box n=\pm\Lambda|u|^2\\
u(0)=u_0,\ n_\pm(0)=n_{\pm0}=n_0\pm i\Lambda^{-1}n_1.
\end{array}\right.\ee
It is obvious that $(u_0,n_{\pm0})\in H^s\times L^2$ by the
regularity of $u_0$, $n_0$ and $n_1$.

Then the expression of energy (or Hamiltonian) above is
\be\label{1.12} H(t)=H(u,n_+)(t)=\|\nabla
u\|_{L^2}^2+\frac{1}{2}\|n_+\|_{L^2}^2+\frac{1}{2}\int(n_++\bar{n}_+)|u|^2dx.\ee

The purpose of us is to imitate the $H^1$ argument with the energy.
But the energy is infinite in the $H^s\times L^2$-setting, thus we
applying a smoothing operator to make $u$ and $n_\pm$ in $H^1\times
H^{1-s}$ and define the usual energy of this new object. However,
the energy is not conserved any more, so the crucial point here is
to estimate the growth of the modified total energy. The main
difficult of this step is the low regularity of $n_+$. In the other
hand the wave equation doesn't possess Strichartz estimates
${\cite{FW}}$ as good as Schr\"{o}dinger, and we need some endpoint
Strichartz estimates, which leads to the requirement of radial
condition.

During the proof, we  find a character of finite time blow up of
Zakharov system, i.e. when $t\rightarrow T^*<\infty$,
$\|Iu(t)\|_{H^1}$ would go to infinite and $\liminf_{t\rightarrow
T^*}\|In_+(t)\|_{H^{1-s}}>0$. In fact, from the local existence
theory, we can get $\|u(t)\|_{H^s}+\|n_+\|_{L^2}\rightarrow\infty$
as $t\rightarrow T^*$, so
$\|Iu(t)\|_{H^1}+\|In_+(t)\|_{H^{1-s}}\rightarrow\infty$. Then we
prove this character by another view of the local existence result
and the particular form of the system, that the nonlinear term of
the second equation is independent on $n_+$.

In Section 2, we'll give some notations, norms and estimates. Then
in Section 3, the local existence theory will be studied while in
Section 4, we'll estimate the change of the modified energy which
is the main part of the paper. In Section 5, the proof for Theorem
\ref{th1} is given.


\section{Notations, Norms and Estimates}

$A\lesssim B$ means there is a universal constant $c>0$, such that
$A\leqslant cB$, and $A\thicksim B$ when both $A\lesssim B$ and
$B\lesssim A$.

$<\xi>=(1+|\xi|^2)^\frac{1}{2}$.

$c+$ means $c+\epsilon$ while $c-$ means $c-\epsilon$, for some
$\epsilon>0$ small enough.

For given $N>>1$, define smoothing operators $I_N:\ $ \be\label{2.1}
\widehat{I_Nf}(\xi)=m_N(\xi)\hat{f}(\xi),\ee where \be\label{2.2}
m_N(\xi)=\left\{
\begin{array}{ll}
1, & |\xi|\leqslant N\\
(\frac{|\xi|}{N})^{s-1}, & |\xi|\geqslant3N,
\end{array}\right.
\ee and $m_N(\xi)$ is smoothing, radial, nonnegative, and monotone
in $|\xi|$. We drop $N$ from the notation for short when there is no
confusion.

By computation, we have \be\label{2.3}
\|f\|_{X^{\tilde{s}_0,b_0}_{\varphi,I}}\lesssim\|If\|_{X^{\tilde{s}_0+1-s,b_0}_{\varphi,I}}\lesssim
N^{1-s}\|f\|_{X^{\tilde{s}_0,b_0}_{\varphi,I}}\ee for any
$\tilde{s}_0\geqslant0$, $b_0\in\mathbb{R}$. Here, we used the
$X^{m,b}_{\varphi,I}$-space which is defined as follows: for an
equation of the form $if_t-\varphi(-i\nabla)f=0$, where $\varphi$ is
a measurable function, let $X^{m,b}_{\varphi}$ be the completion of
$\mathcal{S}(\mathbb{R}\times\mathbb{R}^2)$ with respect to \bea
\|f\|_{X^{m,b}_{\varphi}}:&=&\|<\xi>^m<\tau>^b\mathcal{F}(e^{-it\varphi(-i\partial_x)}f(t,x))\|_{L^2_{\tau\xi}}\nonumber\\
&=&\|<\xi>^m<\tau+\varphi(x)>^b\tilde{f}(\tau,\xi)\|_{L^2_{\tau\xi}}.\label{2.4}\eea
We denote Fourier transform w.r.t both $x$ and $t$ by $\tilde{\ }$,
while only w.r.t $x$ or $t$ by $\hat{\ }$.

For a given time interval $I$, we define
$\|f\|_{X^{m,b}_{\varphi,I}}=\inf_{g|_{I}=f}\|g\|_{X^{m,b}_{\varphi}}$,
and also omit $I$ if there is no confusion.

For $\varphi(\xi)=\pm|\xi|$, we use the notation $X^{m,b}_{\pm}$,
while for $\varphi(\xi)=-|\xi|^2$ simply $X^{m,b}$.

Now, we are listing some well-known estimates for these norms.

1. If $u$ is a solution of $iu_t-\varphi(-i\partial_x)u=0$ with
u(0)=f and $\psi$ is a cut off function in
$C^\infty_0(\mathbb{R})$ with $supp \psi\subset(-2,2)$,
$\psi\equiv1$ on $[-1,1]$, $\psi(t)=\psi(-t)$,
$\psi(t)\geqslant0$, $\psi_\delta(t):=\psi(\frac{t}{\delta})$,
$0<\delta\leqslant1$, we have for $b>0$,
\be\label{2.5}\|\psi_1u\|_{X^{m,b}_{\varphi}}\leqslant
c\|f\|_{H^m}.\ee

If $v$ is a solution of problem $iv_t-\varphi(-i\partial_x)v=F$,
$v(0)=0$, we have for $b'+1\geqslant b\geqslant 0\geqslant
b'>-\frac{1}{2}$ \be\label{2.6} \|\psi_\delta
v\|_{X^{m,b}_{\varphi}}\leqslant
c\delta^{1+b'-b}\|F\|_{X^{m,b'}_{\varphi}}.\ee

The proofs for these two estimates could be found in \cite{GTV}.

2. For $\frac{1}{2}>b>b'\geqslant0$, $0<\delta\leqslant 1$,
$m\in\mathbb{R}$ \be\label{2.7} \|\psi_\delta
f\|_{X^{m,b'}_{\varphi}}\leqslant
c\delta^{b-b'}\|f\|_{X^{m,b}_{\varphi}}.\ee

3. Strichartz estimates.

For $\frac{2}{q}=1-\frac{2}{r}$, $q\geqslant2$ and $u$ is radial,
\be\label{2.8} \|u\|_{L^q_tL^r_x}\leqslant
c\|u\|_{X^{0,\frac{1}{2}+}},\ee and \\
for $\frac{1}{q}=1-\frac{2}{r}$ and $q>2$, \be\label{2.10}
\|v\|_{L^q_tL^r_x}\leqslant c\|v\|_{X^{0,\frac{1}{2}+}_{\pm}}.\ee

4. From \cite{CKSTT1}, one has for $|\xi_i|\sim N_i$, $i=1,2$,
$N_1\leqslant N_2$, \be\label{2.11}
\|u_1u_2\|_{L^2([0,\delta]\times\mathbb{R}^2)}\leqslant
c(\frac{N_1}{N_2})^\frac{1}{2}\|u_1\|_{X^{0,\frac{1}{2}+}}\|u_2\|_{X^{0,\frac{1}{2}+}}.\ee

5. For $s_1\leqslant s_2$
\be\label{2.12}\|f\|_{X^{s_1,b}_\varphi}\leqslant
c\|f\|_{X^{s_2,b}_\varphi},\ee and\\
for $b_1\leqslant b_2$
\be\label{2.13}\|f\|_{X^{s,b_1}_\varphi}\leqslant
c\|f\|_{X^{s,b_2}_\varphi}.\ee

Finally, we give the sharp Gagliardo-Nirenberg inequality for
$\mathbb{R}^2$, which could been found in \cite{W}. \be\label{2.14}
\frac{1}{2}\|u\|_{L^4}^4\leqslant\frac{\|u\|_{L^2}^2}{\|Q\|_{L^2}^2}\|\nabla
u\|_{L^2}^2,\ \ for\ u\in H^1,\ and\ u\neq0.\ee


\section{Local Existence Theory}
The existence and uniqueness for system (\ref{1.11}) holds by the
results of \cite{GTV} for $(u_0,n_{\pm0})\in H^s\times L^2$,
$s\geqslant\frac{1}{2}$.

If we apply operator $I$ to the system (\ref{1.11}), we have
\be\label{3.1} \left\{
\begin{array}{ll}
i\partial_t(Iu)+\triangle Iu=\frac{1}{2}I((n_++n_-)u)\\
(i\partial_t\mp\Lambda)In_\pm=\pm\Lambda I(|u|^2)\\
Iu(0)=Iu_0,\ In_\pm(0)=In_{\pm0}.
\end{array}
\right. \ee

\begin{prop}\label{prop3.1}
Assume $(u_0,n_{\pm0})\in H^s\times L^2$, and $1> s\geqslant
\frac{1}{2}$. Then there exists a positive number
$\delta=min\{(\frac{c}{\|In_{+0}\|_{H^{1-s}}})^{2+17\epsilon},(c\frac{\|In_{+0}\|_{H^{1-s}}}{\|Iu_0\|_{H^1}^2})^{2+17\epsilon}\}$,
with that $\epsilon>0$ is a small enough parameter, such that
system (\ref{3.1}) has a unique local solution ($Iu$, $In_\pm$) in
the time interval $[0,\delta]$ with the
property:\be\label{3.9}\|Iu\|_{X^{1,\frac{1}{2}+}}\lesssim\|Iu_0\|_{H^1},\
\|In_\pm\|_{X^{1-s,\frac{1}{2}+}_\pm}\lesssim\|In_{+0}\|_{H^{1-s}}.\ee
\end{prop}

\begin{proof} Let
\[E=\{(Iu,In_\pm)|\|Iu\|_{X^{1,\frac{1}{2}+}}\lesssim\|Iu_0\|_{H^1},\ \|In_\pm\|_{X^{1-s,\frac{1}{2}+}_\pm}\lesssim\|In_{+0}\|_{H^{1-s}}\},\]
and $(S_0,S_1,S_1)$ defined on $E$ as
\begin{eqnarray*}
&&S_0(Iu)=\psi_1e^{it\triangle}Iu_0-\frac{i}{2}\psi_1\int_0^te^{i(t-s)\triangle}\psi_\delta
I((n_++n_-)u)ds,\\
&&S_1(In_\pm)=\psi_1e^{\mp it\Lambda}In_{\pm0}\mp
i\psi_1\int_0^te^{\mp i(t-s)\Lambda}\psi_\delta\Lambda I(|u|^2)ds,
\end{eqnarray*}
where $\psi_1$ and $\psi_\delta$ are defined before for
$0<\delta\leqslant1$.

Then, taking $b'=-\frac{1}{2}+$ and $b=\frac{1}{2}+$ in (\ref{2.5})
and (\ref{2.6}), it exists \be\label{3.2}
\|S_0(Iu)\|_{X^{1,\frac{1}{2}+}}\leqslant
c\|Iu_0\|_{H^1}+c\|I((n_++n_-)u)\|_{X^{1,-\frac{1}{2}+}},\ee and
\be\label{3.3} \|S_1(In_\pm)\|_{X^{1-s,\frac{1}{2}+}_\pm}\leqslant
c\|In_{\pm0}\|_{H^{1-s}}+c\|\Lambda
I(|u|^2)\|_{X^{1-s,-\frac{1}{2}+}_\pm}.\ee

Next, we use Lemma 3.4 of \cite{GTV} and \cite{CKSTT2} to get
\be\label{3.4}\|I((n_++n_-)u)\|_{X^{1,-\frac{1}{2}+}}\leqslant
c\delta^{\frac{1}{2}-4\epsilon}\|In_\pm\|_{X^{1-s,\frac{1}{2}+}_\pm}\|Iu\|_{X^{1,\frac{1}{2}+}}.\ee

We also use Lemma 3.5 of \cite{GTV} and \cite{CKSTT2} to get
\be\label{3.5}\|\Lambda
I(|u|^2)\|_{X^{1-s,-\frac{1}{2}+}_\pm}\leqslant
c\delta^{\frac{1}{2}-4\epsilon}\|Iu\|_{X^{1,\frac{1}{2}+}}^2.\ee

Combining these estimates together and because $n_+=\bar{n}_-$,
there exists \be\label{3.6}\|S_0(Iu)\|_{X^{1,\frac{1}{2}+}}\leqslant
c\|Iu_0\|_{H^1}+c\delta^{\frac{1}{2}-4\epsilon}\|In_+\|_{X^{1-s,\frac{1}{2}+}_+}\|Iu\|_{X^{1,\frac{1}{2}+}},\ee
and
\be\label{3.7}\|S_1(In_\pm)\|_{X^{1-s,\frac{1}{2}+}_\pm}\leqslant
c\|In_{\pm0}\|_{H^{1-s}}+c\delta^{\frac{1}{2}-4\epsilon}\|Iu\|_{X^{1,\frac{1}{2}+}}^2.\ee

Letting $\delta\sim
\min\{(\frac{1}{\|In_{+0}\|_{H^{1-s+}}})^{2+17\epsilon},(\frac{\|In_{+0}\|_{H^{1-s+}}}{\|Iu_0\|_{H^1}^2})^{2+17\epsilon}\}$,
such that
$\delta^{\frac{1}{2}-4\epsilon}\|In_{+0}\|_{H^{1-s}}\lesssim1$ and
$\delta^{\frac{1}{2}-4\epsilon}\|Iu_0\|_{H^1}^2\lesssim\|In_{+0}\|_{H^{1-s}}$,
then we have\be\label{3.8} (S_0(Iu),S_1(In_\pm))\in E,\ee hence
$(S_0,S_1,S_1):\ E\rightarrow E$.

One can prove $(S_0,S_1,S_1)$ is a contraction map with the same
method. Thus by the standard fixed point theory, we get the local
existence of (\ref{3.1}). And the uniqueness follows in the same
way.
\end{proof}


\begin{prop}\label{prop3.2}
Assume $(u_0,n_{\pm0})\in H^s\times L^2$, with
$1>s\geqslant\frac{1}{2}$, then there exists a positive number
$\tilde{\delta}=\frac{c}{M^{2+17\epsilon}}$, with
$M=(\|u_0\|_{H^s}+\|n_{\pm0}\|_{L^2})$, such that system
(\ref{1.11}) has a unique local solution in the time interval
$[0,\tilde{\delta}]$ with the property: \be\label{3.10}
\|u\|_{X^{s,\frac{1}{2}+}}+\|n_\pm\|_{X^{0,\frac{1}{2}+}_\pm}\lesssim
M.\ee
\end{prop}
\begin{rem}\label{rem3.1}
The proof is almost the same as Proposition \ref{prop3.1}, except
some small changes.
\end{rem}
\begin{proof}
Let $E=\{(u,n_\pm)|\
\|u\|_{X^{s,\frac{1}{2}+}}+\|n_\pm\|_{X^{0,\frac{1}{2}+}_\pm}\lesssim
M\}$ and also define $(S_0,S_1,S_1)$ as
\begin{eqnarray*}
S_0(u)=\psi_1e^{it\triangle}u_0-\frac{i}{2}\psi_1\int_0^te^{i(t-s)\triangle}\psi_{\tilde{\delta}}(n_++n_-)uds,\\
S_1(n_\pm)=\psi_1e^{\mp it\Lambda}n_{\pm0}\mp i\psi_1\int_0^te^{\mp
i(t-s)\Lambda}\psi_{\tilde{\delta}}\Lambda(|u|^2)ds,
\end{eqnarray*}
where $\psi_1$ and $\psi_{\tilde{\delta}}$ are defined as before.

Then, like Proposition \ref{prop3.1} we have \be\label{3.11}
\|S_0(u)\|_{X^{s,\frac{1}{2}+}}\leqslant
c\|u_0\|_{H^s}+c\tilde{\delta}^{\frac{1}{2}-4\epsilon}\|n_+\|_{X^{0,\frac{1}{2}+}_+}\|u\|_{X^{s,\frac{1}{2}+}},\ee
and \be\label{3.12} \|S_1(n_\pm)\|_{X^{0,\frac{1}{2}+}_\pm}\leqslant
c\|n_{\pm0}\|_{L^2}+c\tilde{\delta}^{\frac{1}{2}-4\epsilon}\|u\|_{X^{s,\frac{1}{2}+}}^2.
\ee Thus we just need to take
$\tilde{\delta}=\frac{c}{M^{2+17\epsilon}}$, then the result of the
proposition follows.
\end{proof}


From the above Proposition we can see that,
\begin{cor}\label{cor3.1}
If $(u(t),n_\pm(t))$ is a finite time blow up solution in
$H^s\times L^2$, $\frac{1}{2}\leqslant s<1$, with the initial data
as above,
 then
$\|u(t)\|_{H^s}+\|n_\pm(t)\|_{L^2}\rightarrow\infty$, as
$t\rightarrow T^*$ where $[0,T^*)$ is the maximum life span, which
is also equivalent to
$\|u(t)\|_{H^s}+\|n_+(t)\|_{L^2}\rightarrow\infty$ as
$t\rightarrow T^*$ because of $n_+=\bar{n}_-$.
\end{cor}


\begin{cor}\label{cor3.2}
If $(u(t),n_\pm(t))$ is a finite time blow up solution in
$H^s\times L^2$, $\frac{1}{2}\leqslant s<1$, then
\be\label{3.13}\|Iu(t)\|_{H^1}\rightarrow\infty,\ \ as\
t\rightarrow T^*,\ee and \be\label{3.14} \liminf_{t\rightarrow
T^*}\|In_+(t)\|_{H^{1-s}}>0,\ee i.e. there is a $c>0$ such that
$\|In_+(t)\|_{H^{1-s}}\geqslant c$.
\end{cor}
\begin{proof}
As $\|Iu\|_{H^1}\gtrsim\|u\|_{H^s}$ and
$\|In_+\|_{H^{1-s}}\gtrsim\|n_+\|_{L^2}$, by Corollary \ref{cor3.1},
we have
\be\label{3.18}\|Iu\|_{H^1}+\|In_+\|_{H^{1-s}}\rightarrow\infty,\ \
as\ t\rightarrow T^*,\ for\ fixed\ N>>1.\ee

On the other hand, from the proof of Proposition \ref{prop3.1},
one can find that if replacing $\psi_1$ with $\psi_{T^*}$, the
estimates also hold. Thus, for $T<T^*$, \be\label{3.15}
\|In_+\|_{X^{1-s,\frac{1}{2}+}_+}\leqslant
c\|In_{+0}\|_{H^{1-s}}+cT^{\frac{1}{2}-4\epsilon}\|Iu\|_{X^{1,\frac{1}{2}+}}^2,\ee
and \be\label{3.16} \|Iu\|_{X^{1,\frac{1}{2}+}}\leqslant
c\|Iu_0\|_{H^1}+cT^{\frac{1}{2}-4\epsilon}\|In_+\|_{X^{1-s,\frac{1}{2}+}_+}\|Iu\|_{X^{1,\frac{1}{2}+}}.\ee

Hence, if $\|Iu(t)\|_{H^1}\nrightarrow\infty$, as $t\rightarrow
T^*$, i.e. $\|Iu\|_{L^\infty{[0,T^*),H^1}}\leqslant A$, for some
$A<\infty$, then it has \bea
\|In_+\|_{H^{1-s}}&\lesssim&\|In_+\|_{X^{1-s,\frac{1}{2}+}_+}\lesssim
\|In_{+0}\|_{H^{1-s}}+T^{\frac{1}{2}-4\epsilon}\|Iu\|_{X^{1,\frac{1}{2}+}}^2\nonumber\\
&\lesssim&N^{1-s}\|n_{+0}\|_{L^2}+{T^*}^{\frac{1}{2}-4\epsilon}\|Iu\|_{L^\infty([0,T^*),H^1)}^2\nonumber\\
&\lesssim&N^{1-s}\|n_{+0}\|_{L^2}+{T^*}^{\frac{1}{2}-4\epsilon}A^2<\infty\label{3.17},\eea
for fixed $N>>1$, by Proposition \ref{prop3.1}, (\ref{3.15}) and
$T^*<\infty$, which contradicts to (\ref{3.18}). This proves
(\ref{3.13}).

Next, if $\liminf_{t\rightarrow T^*}\|In_+(t)\|_{H^{1-s}}=0$, then
there would be a subsequence $\{t_n\}$, $t_n\rightarrow T^*$ as
$n\rightarrow\infty$, such that $\lim_{n\rightarrow
\infty}\|In_+(t_n)\|_{H^{1-s}}=0$, so from (\ref{3.16}) we have,
\begin{eqnarray*}\|Iu(t_n)\|_{X^{1,\frac{1}{2}+}}&\leqslant&
c\|Iu_0\|_{H^1}+c{t_n}^{\frac{1}{2}-4\epsilon}\|In_+(t_n)\|_{X^{1-s,\frac{1}{2}+}_+}\|Iu(t_n)\|_{X^{1,\frac{1}{2}+}}\\
&\leqslant&cN^{1-s}\|u_0\|_{H^s}+c{T^*}^{\frac{1}{2}-4\epsilon}\|In_+(\tilde{t}_n)\|_{H^{1-s}}\|Iu(t_n)\|_{X^{1,\frac{1}{2}+}},
\end{eqnarray*} for some $\tilde{t}_n,$ which satisfies $|\tilde{t}_n-t_n|\lesssim\delta$
by the local existence theory Proposition \ref{prop3.1}. Hence,
since $T^*<\infty$, for $n\rightarrow\infty$, \be\label{3.19}
\|Iu(t_n)\|_{X^{1,\frac{1}{2}+}}\lesssim N^{1-s}\|u_0\|_{H^s}.\ee
(\ref{3.19}) gives
\[\|Iu(t_n)\|_{H^1}<c<\infty,\] for fixed $N>>1$, which contradicts to (\ref{3.13}).

\end{proof}


\section{Estimates for the Modified Energy}
In this section we'll get the exact control of the increment of the
modified energy.

As the modified energy is $H(t)=H(Iu,In_+)=\|\nabla
Iu\|_{L^2}^2+\frac{1}{2}\|In_+\|_{L^2}^2+\frac{1}{2}\int
I(n_++\bar{n}_+)|Iu|^2dx$, and it is not conserved any more, we have
to control its growth. The following is the main proposition of the
paper:

\begin{prop}\label{prop4.1}
Let $(Iu,In_\pm)$ be a solution of (\ref{3.1}) on $[0,\delta]$ in
the sense of Proposition \ref{prop3.1}. Then the following estimate
holds ($N>>1$):\be\label{4.2} |H(\delta)-H(0)|\leqslant
                    cN^{-2+s+}\delta^{0+}\|In_+\|_{X^{1-s,\frac{1}{2}+}_+}\|Iu\|_{X^{1,\frac{1}{2}+}}^2
                    +cN^{-2+s+}\|In_+\|_{X^{1-s,\frac{1}{2}+}_+}^2\|Iu\|_{X^{1,\frac{1}{2}+}}^2.\ee
\end{prop}

\begin{proof}\begin{eqnarray*}\frac{dH(t)}{dt}=&&2Re\int\nabla\overline{Iu}\nabla
Iu_tdx+Re\int
\overline{In}_+In_{+t}+\frac{1}{2}\int((In_+)_t+(In_-)_t)|Iu|^2+Re\int(In_++In_-)\overline{Iu}Iu_t\\
=&&-Im\int\triangle \overline{Iu}(I((n_++n_-)u)-(In_++In_-)Iu)\\
&&+\frac{1}{2}Im\int \overline{I((n_++n_-)u)}(I((n_++n_-)u)-(In_++In_-)Iu)\\
&&-Im\int \overline{In_+}\Lambda(I(|u|^2)-|Iu|^2).\\
 \end{eqnarray*}
Integrate by $t$ on $[0,\delta)$, it has \bea
|H(\delta)-H(0)|\leqslant&&|\int_0^\delta\int_{\mathbb{R}^2}\triangle
\overline{Iu}(I((n_++n_-)u)-(In_++In_-)Iu)dxdt|\nonumber\\
&&+\frac{1}{2}|\int_0^\delta\int_{\mathbb{R}^2}\overline{I((n_++n_-)u)}(I((n_++n_-)u)-(In_++In_-)Iu)dxdt|\nonumber\\
&&+|\int_0^\delta\int_{\mathbb{R}^2}\overline{In_+}\Lambda(I(|u|^2)-|Iu|^2)dxdt|\nonumber\\
=&&I+II+III.\label{4.1}\eea

To prove Proposition \ref{prop4.1}, we have to control $I$, $II$ and
$III$ in (\ref{4.1}) respectively.


First for $I$, it has,
\begin{lem}\label{lem1}
$I\lesssim(N^{-2+s+}\delta^{0+}+N^{-\frac{5}{2}+s+2\epsilon}\delta^{\frac{1}{2}-})
\|In_+\|_{X^{1-s+,\frac{1}{2}+}_+}\|Iu\|_{X^{1,\frac{1}{2}}}^2$.
\end{lem}
\begin{proof}
As \begin{eqnarray*} I&=&|\int_0^\delta\int_{\mathbb{R}^2}\triangle
\overline{Iu}(I((n_++n_-)u)-(In_++In_-)Iu)dxdt|\\
                     &\sim&|\int_0^\delta\int_*|\xi_1|^2
m(\xi_1)\hat{\bar{u}}(\xi_1)(\frac{m(\xi_2+\xi_3)-m(\xi_2)m(\xi_3)}{m(\xi_2)m(\xi_3)})m(\xi_2)\hat{n}_+(\xi_2)m(\xi_3)\hat{u}(\xi_3))d\xi
dt|,
\end{eqnarray*}
here $*$ denotes integration over the set $\{\sum_{i=1}^3\xi_i=0\}$
(or $\{\sum_{i=1}^4\xi_i=0\}$).

We break the function $u$ and $n_+$ into a sum of dyadic
constituents, each with frequency support $<\xi_i>\sim 2^j$,
$j=0,\cdots$ and denote $u_i=P_{N_i}Iu$, $n_{+i}=P_{N_i}In_+$.

In the following, let's note $m_i=m(\xi_i)$, $|\xi_i|=N_i$,
$N_{max}=\max_{1\leqslant i\leqslant3}N_i$ (or
$N_{max}=\max_{1\leqslant i\leqslant4}N_i$).

Remark also that w.l.o.g. $\hat{u}_1$, $\hat{n}_{+2}$,
$\hat{u}_3\geqslant0$.

As if both $N_2$ and $N_3<<N$, then $m_2=m_3=1$ such that
$\frac{m(\xi_2+\xi_3)-m_2m_3}{m_2m_3}=0$, the left hand side of the
inequality becomes 0, which is trivial.

Inclusively, we just need to prove \be\label{4.3}
I':=N_1^2\int_0^\delta\int_*
|\frac{m(\xi_2+\xi_3)-m_2m_3}{m_2m_3}|\hat{u}_1\hat{n}_{+2}\hat{u}_3\lesssim
N^{-2+s+}\delta^{0+}\|In_+\|_{X^{1-s,\frac{1}{2}+}_+}\|Iu\|^2_{X^{1,\frac{1}{2}+}},
\ee with the assumption that at least one of $N_2$ and $N_3\gtrsim
N$.


{\bf Case 1.} $N_{max}\sim N_2\sim N_3\gtrsim N$

In this case,
$\frac{m(\xi_2+\xi_3)-m_2m_3}{m_2m_3}\lesssim\frac{m_1}{m_2m_3}\lesssim\frac{1}{m_2m_3}\sim
\frac{1}{m_2^2}\lesssim(\frac{N_2}{N})^{2(1-s)}$, by
$\sum_{i=1}^3\xi_i=0$ and the definition of $m$.


{\bf a.} $N_1\leqslant N_3^\epsilon$.

It exists, \bea I'&\lesssim&
N_1^2(\frac{N_2}{N})^{2(1-s)}\|n_{+2}\|_{L^2_{t,x}}\|u_1u_3\|_{L^2_{t,x}}\nonumber\\
&\lesssim&
N_1^2(\frac{N_2}{N})^{2(1-s)}\delta^{\frac{1}{2}-}\|n_{+2}\|_{X^{0,\frac{1}{2}+}_+}(\frac{N_1}{N_3})^\frac{1}{2}\|u_1\|_{X^{0,\frac{1}{2}+}}\|u_3\|_{X^{0,\frac{1}{2}+}}\nonumber\\
&\lesssim&N_1^2(\frac{N_2}{N})^{2(1-s)}\delta^{\frac{1}{2}-}(\frac{N_1}{N_3})^\frac{1}{2}\frac{1}{<N_1>}\frac{1}{N_2^{1-s}}\frac{1}{N_3}\|n_{+2}\|_{X^{1-s,\frac{1}{2}+}_+}\|u_1\|_{X^{1,\frac{1}{2}+}}\|u_3\|_{X^{1,\frac{1}{2}+}}\nonumber\\
&\lesssim&N_{max}^{-\epsilon}N^{-\frac{5}{2}+s+2\epsilon}\delta^{\frac{1}{2}-}\|In_+\|_{X^{1-s,\frac{1}{2}+}_+}\|Iu\|_{X^{1,\frac{1}{2}+}}^2,\label{4.4}
\eea by the definition of $X^{m,b}_\varphi$-space, (\ref{2.7}),
(\ref{2.11}) and (\ref{2.12}).


{\bf b.} $N_2\sim N_1\sim N_3$

We have to take the fourier transform of $t$ into account in this
case, and w.l.o.g $\tilde{u}_1$, $\tilde{n}_{+2}$, $\tilde{u}_3>0$.

\be\label{4.5} I'\lesssim
N_1^2(\frac{N_2}{N})^{2(1-s)}\int_{**}\tilde{u}_1(\tau_1)\tilde{n}_{+2}(\tau_2)\tilde{u}_3(\tau_3)|\hat{\phi}(\tau_0)|d\xi
d\tau,\ee here $**$ denotes integration over
$\sum_{i=1}^3\xi_i=\sum_{i=0}^3\tau_i=0$, and $\phi(t)$ is the
characteristic function of the time interval $[0,\delta]$.

It is known that
$\hat{\phi}(\tau)=\frac{1}{\sqrt{2\pi}}\frac{e^{i\tau\delta}-1}{i\tau}\in
L_\tau^{1+}$ but not in $L_\tau^1$.

To deal with this case, we need the following algebraic inequality.
\be\label{4.6}|\xi_1|\lesssim<\tau_1+|\xi_1|^2>^\frac{1}{2}+<\tau_2+|\xi_2|>^\frac{1}{2}+<\tau_3+|\xi_3|^2>^\frac{1}{2}+|\tau_0|^\frac{1}{2},\ee
and consider every 4 cases according to which terms on the r.h.s is
dominant.


{\bf Subcase 1.} $<\tau_1+|\xi_1|^2>^\frac{1}{2}$ dominant.

\bea
I'&\lesssim&(\frac{N_2}{N})^{2(1-s)}\int_{**}<\tau_1+|\xi_1|^2>^{\frac{1}{2}+}|\xi_1|^{1-}\tilde{u}_1\tilde{n}_{+2}\tilde{u}_3|\hat{\phi}|\nonumber\\
  &\lesssim&(\frac{N_2}{N})^{2(1-s)}N_1^{1-}\|u_1\|_{X^{0,\frac{1}{2}+}}\|\mathcal{F}^{-1}(|\hat\phi|)n_+u_3\|_{L^2_{t,x}}\nonumber\\
  &\lesssim&(\frac{N_2}{N})^{2(1-s)}\frac{1}{N_1^\epsilon}\|u_1\|_{X^{1,\frac{1}{2}+}}\|\mathcal{F}^{-1}(|\hat{\phi}|)\|_{L_t^{\infty-}}\|n_{+2}\|_{L^\infty_tL^2_x}\|u_3\|_{L_t^{2+}L^\infty_x}\nonumber\\
  &\lesssim&(\frac{N_2}{N})^{2(1-s)}\frac{1}{N_1^\epsilon}\delta^{0+}\|u_1\|_{X^{1,\frac{1}{2}+}}\|n_{+2}\|_{X^{0,\frac{1}{2}+}_+}N_3^\frac{\epsilon}{2}\|u_3\|_{X^{0,\frac{1}{2}+}}\nonumber\\
  &\lesssim&(\frac{N_2}{N})^{2(1-s)}\frac{1}{N_1^\epsilon}\delta^{0+}\frac{1}{N_2^{1-s}}\frac{1}{N_3}N_3^\frac{\epsilon}{2}\|In_+\|_{X^{1-s,\frac{1}{2}+}_+}\|Iu\|_{X^{1,\frac{1}{2}+}}^2\nonumber\\
  &\lesssim&N_{max}^{-\frac{\epsilon}{2}}N^{-2+s}\delta^{0+}\|In_+\|_{X^{1-s,\frac{1}{2}+}_+}\|Iu\|_{X^{1,\frac{1}{2}+}}^2,\label{4.7}
\eea by H\"{o}lder inequality, Berstein inequality, (\ref{2.8}),
(\ref{2.10}), and Hausdorff-Young, which gives
\be\label{4.8}\|\mathcal{F}^{-1}(|\hat{\phi}|)\|_{L_t^{\infty-}}\lesssim\|\hat{\phi}\|_{L^{1+}_\tau}\lesssim\delta^{0+}.
\ee


{\bf Subcase 2.} $<\tau_2+|\xi_2|>^\frac{1}{2}$ dominant.

\bea
I'&\lesssim&(\frac{N_2}{N})^{2(1-s)}N_1^{1-}\int_{**}<\tau_2+|\xi_2|>^{\frac{1}{2}+}|\hat{\phi}(\tau_0)|\tilde{u}_1(\tau_1)\tilde{n}_{+2}(\tau_2)\tilde{u}_3(\tau_3)d\xi
d\tau\nonumber\\
  &\lesssim&(\frac{N_2}{N})^{2(1-s)}N_1^{1-}\|n_{+2}\|_{X^{0,\frac{1}{2}+}_+}\|\mathcal{F}^{-1}(|\hat{\phi}|)u_1u_3\|_{L^2_{t,x}}\nonumber\\
  &\lesssim&(\frac{N_2}{N})^{2(1-s)}N_1^{1-}\frac{1}{N_2^{1-s}}\|n_{+2}\|_{X^{1-s,\frac{1}{2}+}_+}\|\mathcal{F}^{-1}(\hat{\phi})\|_{L_t^{\infty-}}\|u_1u_3\|_{L^{2+}_{t,x}}\nonumber\\
  &\lesssim&(\frac{N_2}{N})^{2(1-s)}N_1^{1-}\frac{1}{N_2^{1-s}}\|n_{+2}\|_{X^{1-s,\frac{1}{2}+}_+}\delta^{0+}\frac{1}{N_1}\frac{1}{N_3}N_3^{\frac{\epsilon}{2}}\|u_1\|_{X^{1,\frac{1}{2}+}}\|u_3\|_{X^{1,\frac{1}{2}}}\nonumber\\
  &\lesssim&(\frac{N_2}{N})^{2(1-s)}N_1^{1-}\frac{1}{N_2^{1-s}}\delta^{0+}\frac{1}{N_1}\frac{1}{N_3}N_3^{\frac{\epsilon}{2}}\|In_+\|_{X^{1-s,\frac{1}{2}+}_+}\|Iu\|_{X^{1,\frac{1}{2}+}}\nonumber\\
  &\lesssim&N_{max}^{-\frac{\epsilon}{2}}N^{-2+s}\delta^{0+}\|In_+\|_{X^{1-s,\frac{1}{2}+}_+}\|Iu\|_{X^{1,\frac{1}{2}+}},\label{4.9}
\eea for the same reason as subcase 1.


{\bf Subcase 3.} $<\tau_3+|\xi|^2>^\frac{1}{2}$ dominant.

Almost the same as subcase 1.


{\bf subcase 4.} $|\tau_0|^\frac{1}{2}$ dominant.

\bea
I'&\lesssim&(\frac{N_2}{N})^{2(1-s)}N_1\int_{**}|\tau_0|^{\frac{1}{2}}|\hat{\phi}(\tau_0)|\tilde{u}_1\tilde{n}_{+2}\tilde{u}_3d\xi
d\tau\nonumber\\
  &\lesssim&(\frac{N_2}{N})^{2(1-s)}N_1\|\tilde{u}_1\|_{L^2_{\xi_1}L^{1+}_{\tau_1}}\||\tau|^\frac{1}{2}|\hat{\phi}|\ast\tilde{n}_{+2}\ast\tilde{u}_3\|_{L^2_{\xi}L^{1-}_\tau}.\label{4.10}
\eea

The first factor is estimated as follows by H\"{o}lder w.r.t.
$\tau_1$: \bea
\|\tilde{u}_1\|_{L^2_{\xi_1}L^{1+}_{\tau_1}}&=&\|\tilde{u}_1<\tau_1+|\xi_1|^2>^{\frac{1}{2}+}<\tau_1+|\xi_1|^2>^{-\frac{1}{2}-}\|_{L^2_{\xi_1}L^{1+}_{\tau_1}}\nonumber\\
                                            &\lesssim&\|\tilde{u}_1<\tau_1+|\xi_1|^2>^{\frac{1}{2}+}\|_{L^2_{\xi_1}L^2_{\tau_1}}\|<\tau_1+|\xi_1|^2>^{-\frac{1}{2}-}\|_{L^\infty_{\xi_1}L^{\frac{2(1+\epsilon)}{1-\epsilon}}_{\tau_1}}\nonumber\\
                                            &\lesssim&\frac{1}{N_1}\|u_1\|_{X^{1,\frac{1}{2}+}},\label{4.11}
\eea since
$(-\frac{1}{2}-\epsilon)\frac{2(1+\epsilon)}{1-\epsilon}<-1$ which
ensure the integrable condition at infinite for $\tau_1$.

The second factor is bounded by Young'inequality by

\bea
\||\hat{\phi}|\ast\tilde{n}_{+2}\ast\tilde{u}_3\|_{L^2_{\xi}L^{1-}_\tau}&\lesssim&\||\tau|^\frac{1}{2}|\hat{\phi}|\|_{L_\tau^{2+\epsilon}}\|\tilde{u}_3\ast\tilde{n}_{+2}\|_{L_\xi^2L_\tau^{2-2\epsilon}}\nonumber\\
&\lesssim&\delta^{0+}\|\tilde{u}_3\|_{L_\xi^1L_\tau^{2-2\epsilon}}\|\tilde{n}_{+2}\|_{L_\xi^2L_\tau^1},\label{4.12}
\eea here we use the bound
$\||\tau|^\frac{1}{2}|\hat{\phi}|\|_{L^{2+\epsilon}_{\tau}}\lesssim\delta^{0+}$.

Because
\bea\|\tilde{u}_3\|_{L_\xi^1L_\tau^{2-2\epsilon}}&\lesssim&N_3^\frac{\epsilon}{2}\|\tilde{u}_3<\xi_3><\tau_3+|\xi_3|^2>^{\frac{1}{2}+}<\xi_3>^{-1-\frac{\epsilon}{2}}<\tau_3+|\xi_3|^2>^{-\frac{1}{2}-}\|_{L_\xi^2L_\tau^{2-2\epsilon}}\nonumber\\
                                                 &\lesssim&N_3^\frac{\epsilon}{2}\|\tilde{u}_3<\xi_3><\tau_3+|\xi_3|^2>^{\frac{1}{2}+}\|_{L_\xi^2L_\tau^2}\|<\xi_3>^{-1-\frac{\epsilon}{2}}<\tau_3+|\xi_3|^2>^{-\frac{1}{2}-}\|_{L_\xi^2L_\tau^{\frac{2(1-\epsilon)}{\epsilon}}}\nonumber\\
                                                 &\lesssim&N_3^\frac{\epsilon}{2}\|u_3\|_{X^{0,\frac{1}{2}+}}\lesssim
                                                 N_3^{\frac{\epsilon}{2}-1}\|u_3\|_{X^{1,\frac{1}{2}+}},\label{4.13}
\eea by $\frac{2(1-\epsilon)}{\epsilon}(-\frac{1}{2}-\epsilon)<-1$
and $2(-1-\frac{\epsilon}{2})+1<-1$.

And \bea
\|\tilde{n}_{+2}\|_{L_\xi^2L_\tau^1}&=&\|\tilde{n}_{+2}<\tau_2+|\xi_2|>^{\frac{1}{2}+}<\tau_2+|\xi_2|>^{-\frac{1}{2}-}\|_{L_\xi^2L_\tau^1}\nonumber\\
&\lesssim&\|\tilde{n}_{+2}<\tau_2+|\xi_2|>^{\frac{1}{2}+}\|_{L_\xi^2L_\tau^2}\|<\tau_2+|\xi_2|>^{-\frac{1}{2}-}\|_{L_\xi^\infty
L_\tau^2}\nonumber\\
&\lesssim&\|n_{+2}\|_{X^{0,\frac{1}{2}+}_+}\lesssim
\frac{1}{N_2^{1-s}}\|n_{+2}\|_{X^{1-s,\frac{1}{2}+}_+},\label{4.14}
\eea since $2(-\frac{1}{2}-\epsilon)<-1$.

Hence, \be\label{4.15} (\ref{4.12})\lesssim
N_3^{\frac{\epsilon}{2}-1}N_2^{-(1-s)}\delta^{0+}\|u_3\|_{X^{1,\frac{1}{2}+}}\|n_{+2}\|_{X^{1-s,\frac{1}{2}+}_+},
\ee then with (\ref{4.11}), it has \be\label{4.16} I'\lesssim
(\frac{N_2}{N})^{1-s}N_3^{\frac{\epsilon}{2}-1}\frac{1}{N_2^{1-s}}\delta^{0+}\|In_+\|_{X^{1-s,\frac{1}{2}+}_+}\|Iu\|_{X^{1,\frac{1}{2}+}}^2\lesssim
N_{max}^{-\frac{\epsilon}{2}}N^{-2+s+\epsilon}\delta^{0+}\|In_+\|_{X^{1-s,\frac{1}{2}+}_+}\|Iu\|_{X^{1,\frac{1}{2}+}}^2.
\ee

{\bf c.} $N_3^\epsilon\leqslant N_1\leqslant N_3\sim N_2$.

Deal with this situation like case b, hence,

for subcase1, \bea
I'&\lesssim&(\frac{N_2}{N})^{2(1-s)}\frac{1}{N_1^\epsilon}\frac{1}{N_2^{1-s}}\frac{1}{N_3}N_3^\frac{\epsilon}{2}\delta^{0+}\|In_+\|_{X^{1-s,\frac{1}{2}+}}\|Iu\|_{X^{1,\frac{1}{2}+}}^2\nonumber\\
&\lesssim&N_{max}^{-\frac{\epsilon}{3}}N^{-2+s+\epsilon}N_1^{-\epsilon}N_3^{-\frac{\epsilon}{6}}\delta^{0+}\|In_+\|_{X^{1-s,\frac{1}{2}+}}\|Iu\|_{X^{1,\frac{1}{2}+}}^2\nonumber\\
&\lesssim&N_{max}^{-\frac{\epsilon}{3}}N^{-2+s+\epsilon}N_3^{-\epsilon^2}N_3^{-\frac{\epsilon}{6}}\delta^{0+}\|In_+\|_{X^{1-s,\frac{1}{2}+}}\|Iu\|_{X^{1,\frac{1}{2}+}}^2\nonumber\\
&\lesssim&N_{max}^{-\frac{\epsilon}{3}}N^{-2+s+\epsilon}N_3^{-\epsilon^2-\frac{\epsilon}{6}}\delta^{0+}\|In_+\|_{X^{1-s,\frac{1}{2}+}}\|Iu\|_{X^{1,\frac{1}{2}+}}^2\nonumber\\
&\lesssim&N_{max}^{-\frac{\epsilon}{3}}N^{-2+s+\epsilon}\|In_+\|_{X^{1-s,\frac{1}{2}+}}\|Iu\|_{X^{1,\frac{1}{2}+}}^2.\nonumber
\eea

And the other three subcases could be dealt with in the same way.


{\bf Case 2} $N_{max}\sim N_2\sim N_1\gtrsim N$.

{\bf a} $N_1\sim N_2\sim N_3$.

This case is the same as Case 1b.

{\bf b} $N_1\sim N_2>>N_3\gtrsim N$.

Then
$|\frac{m(\xi_2+\xi_3)-m_2m_3}{m_2m_3}|\lesssim\frac{m(\xi_2+\xi_3)}{m_2m_3}\lesssim\frac{m_2}{m_2m_3}=\frac{1}{m_3}\lesssim(\frac{N_3}{N})^{1-s}$,
and with the same assumption and argument as Case 1b gives
\be\label{4.17} I'\lesssim
(\frac{N_3}{N})^{1-s}N_1^2\int_{**}\tilde{u}_1(\tau_1)\tilde{n}_{+2}(\tau_2)\tilde{u}_3(\tau_3)|\hat{\phi(\tau_0)}|d\xi
d\tau.\ee

We also divide this case into 4 subcases as before.

As the main part is almost the same, we only show some difference in
the follow.

{\bf Subcase 1} $<\tau_1+|\xi_1|^2>^\frac{1}{2}$ dominant.

\bea
I'&\lesssim&(\frac{N_3}{N})^{1-s}N_1^2\int_{**}\tilde{u}_1\tilde{n}_{+2}\tilde{u}_3|\hat{\phi}|d\xi d\tau\nonumber\\
  &\lesssim&(\frac{N_3}{N})^{1-s}N_1^{1-}\|u_1\|_{X^{0,\frac{1}{2}+}}\|n_{+2}u_3\mathcal{F}^{-1}(|\hat{\phi}|)\|_{L^2_{t,x}}\nonumber\\
  &\lesssim&(\frac{N_3}{N})^{1-s}\frac{1}{N_1^\epsilon}\delta^{0+}\frac{1}{N_2^{1-s}}\frac{N_3^\frac{\epsilon}{2}}{N_3}\|In_+\|_{X^{1-s,\frac{1}{2}+}_+}\|Iu\|_{X^{1,\frac{1}{2}+}}^2\nonumber\\
  &\lesssim&N_{max}^{-\frac{\epsilon}{2}}N^{-2+s}\delta^{0+}\|In_+\|_{X^{1-s,\frac{1}{2}+}_+}\|Iu\|_{X^{1,\frac{1}{2}+}}^2.
  \eea


{\bf Subcase 2} $<\tau_2+|\xi_2|>^\frac{1}{2}$ dominant.

\begin{eqnarray*}
I'&\lesssim&(\frac{N_3}{N})^{1-s}\frac{1}{N_1^\epsilon}\delta^{0+}\frac{1}{N_2^{1-s}}\frac{1}{N_3}N_3^\frac{\epsilon}{2}\|In_+\|_{X^{1-s,\frac{1}{2}+}_+}\|Iu\|_{X^{1,\frac{1}{2}+}}^2\\
  &\lesssim&N_{max}^{-\frac{\epsilon}{2}}N^{-2+s}\delta^{0+}\|In_+\|_{X^{1-s,\frac{1}{2}+}_+}\|Iu\|_{X^{1,\frac{1}{2}+}}^2.
\end{eqnarray*}


{\bf Subcase 3} $<\tau_3+|\xi_3|^2>^\frac{1}{2}$ dominant.

\begin{eqnarray*}
I'&\lesssim&(\frac{N_3}{N})^{1-s}\frac{1}{N_1^\epsilon}\delta^{0+}\frac{1}{N_2^{1-s}}\frac{1}{N_3}N_1^\frac{\epsilon}{2}\|In_+\|_{X^{1-s,\frac{1}{2}+}_+}\|Iu\|_{X^{1,\frac{1}{2}+}}^2\\
  &\lesssim&N_{max}^{-\frac{\epsilon}{2}}N^{-2+s}\delta^{0+}\|In_+\|_{X^{1-s,\frac{1}{2}+}_+}\|Iu\|_{X^{1,\frac{1}{2}+}}^2.
\end{eqnarray*}


{\bf Subcase 4} $|\tau_0|^\frac{1}{2}$ dominant.

\begin{eqnarray*}
I'&\lesssim&(\frac{N_3}{N})^{1-s}N_3^{\frac{\epsilon}{2}-1}\frac{1}{N_2^{1-s}}\delta^{0+}\|In_+\|_{X^{1-s,\frac{1}{2}+}_+}\|Iu\|_{X^{1,\frac{1}{2}+}}^2\\
  &\lesssim&N_{max}^{-\epsilon}N^{-2+s+\epsilon}\delta^{0+}\|In_+\|_{X^{1-s,\frac{1}{2}+}_+}\|Iu\|_{X^{1,\frac{1}{2}+}}^2.
\end{eqnarray*}


{\bf c} $N_1\sim N_2\gtrsim N>>N_3$.

$|\frac{m(\xi_2+\xi_3)-m_2m_3}{m_2m_3}|\sim|\frac{m(\xi_2+\xi_3)-m_2}{m_3}|\lesssim|\frac{\nabla
m_2\xi_3}{m_2}|\lesssim|\frac{N_3}{N_2}|$. We also deal with it in
4 subcases.

{\bf Subcase 1} $<|\tau_1+|\xi_1|^2>^\frac{1}{2}$ dominant.
\be\label{4.18}I'\lesssim
\frac{N_3}{N_2}\frac{1}{N_1^\epsilon}N_3^\frac{\epsilon}{2}\frac{1}{<N_3>}\frac{1}{N_2^{1-s}}\delta^{0+}\|In_+\|_{X^{1-s,\frac{1}{2}+}_+}\|Iu\|_{X^{1,\frac{1}{2}+}}^2
\lesssim
N_{max}^{-\frac{\epsilon}{2}}N^{-2+s}\delta^{0+}\|In_+\|_{X^{1-s+,\frac{1}{2}+}_+}\|Iu\|_{X^{1,\frac{1}{2}+}}^2,\ee
by H\"{o}lder inequality, Berstein inequality, (\ref{2.8}),
(\ref{2.10}) and (\ref{4.8}).

{\bf Subcase 2} $<\tau_2+|\xi_2|>^\frac{1}{2}$ dominant.

\be\label{4.19}I'\lesssim
\frac{N_3}{N_2}\frac{1}{N_1^\epsilon}\delta^{0+}\frac{1}{<N_3>}N_3^{\frac{\epsilon}{2}}\frac{1}{N_2^{1-s}}\|In_+\|_{X^{1-s,\frac{1}{2}+}_+}\|Iu\|_{X^{1,\frac{1}{2}+}}^2\lesssim
N_{max}^{-\frac{\epsilon}{2}}N^{-2+s}\delta^{0+}\|In_+\|_{X^{1-s,\frac{1}{2}+}_+}\|Iu\|_{X^{1,\frac{1}{2}+}},\ee
also by H\"{o}lder inequality, Berstein inequality, (\ref{2.8}),
(\ref{2.10}) and (\ref{4.8}).

{\bf Subcase 3} $<\tau_3+|\xi_3|^2>^\frac{1}{2}$ dominant.

The same as subcase 1.

{\bf Subcase 4} $|\tau_0|^\frac{1}{2}$ dominant.

By H\"{o}lder inequality, (\ref{4.11}) and (\ref{4.15}), we have
\be\label{4.20}
I'\lesssim\frac{N_3}{N_2}\frac{<N_3>^\frac{\epsilon}{2}}{<N_3>}\frac{1}{N_2^{1-s}}\delta^{0+}\|In_+\|_{X^{1-s,\frac{1}{2}+}_+}\|Iu\|_{X^{1,\frac{1}{2}+}}^2\lesssim
N_{max}^{-\frac{\epsilon}{2}}N^{-2+s}\delta^{0+}\|In_+\|_{X^{1-s,\frac{1}{2}+}_+}\|Iu\|_{X^{1,\frac{1}{2}+}}^2.
\ee

{\bf Case 3} $N_{max}\sim N_3\sim N_1\gtrsim N$.

{\bf a} $N_3\sim N_1\sim N_3$.

The same as Case 1b.


{\bf b} $N_3\sim N_1>>N_2\gtrsim N$.

$|\frac{m(\xi_2+\xi_3)-m_2m_3}{m_2m_3}|\lesssim\frac{m(\xi_2+\xi_3)}{m_2m_3}\lesssim\frac{m_3}{m_2m_3}=\frac{1}{m_2}\lesssim(\frac{N_2}{N})^{1-s}$.

Argue as before, we have:
\begin{eqnarray*}
I'&\lesssim&(\frac{N_2}{N})^{1-s}(\frac{1}{N_1^\epsilon}\frac{1}{N_2^{1-s}}\frac{1}{N_3}N_3^{\frac{\epsilon}{2}}+\frac{1}{N_1^\epsilon}\frac{1}{N_2^{1-s}}\frac{1}{N_3}N_1^{\frac{\epsilon}{2}}+N_3^{\frac{\epsilon}{2}-1}\frac{1}{N_2^{1-s}})\delta^{0+}\|In_+\|_{X^{1-s,\frac{1}{2}+}_+}\|Iu\|_{X^{1,\frac{1}{2}+}}^2\\
  &\lesssim&N_{max}^{-\frac{\epsilon}{2}}N^{-2+s+}\delta^{0+}\|In_+\|_{X^{1-s,\frac{1}{2}+}_+}\|Iu\|_{X^{1,\frac{1}{2}+}}^2.
\end{eqnarray*}


{\bf c} $N_3\sim N_1\gtrsim N>>N_2$.

Then
$|\frac{m(\xi_2+\xi_3)-m_2m_3}{m_3m_2}|\sim|\frac{m(\xi_2+\xi_3)-m_3}{m_3}|\lesssim|\frac{\nabla
m_3\xi_2}{m_3}|\lesssim\frac{N_2}{N_3}$. And

\begin{eqnarray*}
I'&\lesssim&\frac{N_2}{N_3}(\frac{1}{N_1^\epsilon}\frac{1}{<N_2>^{1-s}}\frac{1}{N_3}N_3^\frac{\epsilon}{2}+\frac{1}{N_1\epsilon}\frac{1}{<N_2>^{1-s}}\frac{1}{N_3}N_1^\frac{\epsilon}{2}+N_3^{\frac{\epsilon}{2}-1}\frac{1}{<N_2>^{1-s}})\delta^{0+}\|In_+\|_{X^{1-s,\frac{1}{2}+}_+}\|Iu\|_{X^{1,\frac{1}{2}+}}^2\\
  &\lesssim&N_{max}^{-\frac{\epsilon}{2}}N^{-2+s+}\delta^{0+}\|In_+\|_{X^{1-s,\frac{1}{2}+}_+}\|Iu\|_{X^{1,\frac{1}{2}+}}^2.
  \end{eqnarray*}

Taking all the above estimates into account, the result of Lemma
\ref{lem1} holds.

\end{proof}


Now, let's deal with $II$.

\begin{lem}\label{lem2}
\be\label{4.23} II\lesssim N^{-2+s+}
\|In_+\|_{X^{1-s,\frac{1}{2}+}_+}^2\|Iu\|_{X^{1,\frac{1}{2}+}}^2.
\ee
\end{lem}
\begin{proof}
Like part $I$, to prove the estimate for $II$, we just need to prove
\be\label{4.24}
II':=\int_0^\delta\int_*|\frac{m(\xi_1+\xi_2)}{m_1m_2}||\frac{m(\xi_3+\xi_4)-m_3m_4}{m_3m_4}|\hat{n}_{+1}\hat{u}_2\hat{n}_{+3}\hat{u}_4\lesssim
N^{-2+s+}\|In\|_{X^{1-s,\frac{1}{2}+}_+}^2\|Iu\|_{X^{1,\frac{1}{2}+}}^2,
\ee with the same notations and assumptions in Lemma \ref{lem1}.

One can also easy to see that if both $N_3$ and $N_4<<N$, l.h.s
would be zero, which is trivial, so we can suppose at least one of
$N_3$ and $N_4\gtrsim N$.

{\bf Case 1} $N_{max}\sim N_3\sim N_4\gtrsim N$.

So
$|\frac{m(\xi_3+\xi_4)-m_3m_4}{m_3m_4}|\lesssim|\frac{m(\xi_3+\xi_4)}{m_3m_4}|\lesssim\frac{1}{m_3m_4}\sim\frac{1}{m_3^2}\lesssim(\frac{N_3}{N})^{2(1-s)}$.

{\bf a} $N_1$, $N_2<<N$.

In this case, $|\frac{m(\xi_1+\xi_2)}{m_1m_2}|\sim1$, and \bea
II'&\lesssim&
(\frac{N_3}{N})^{2(1-s)}\|n_{+1}u_2\|_{L^2_{t,x}}\|n_{+3}u_4\|_{L^2_{t,x}}\nonumber\\
&\lesssim&(\frac{N_3}{N})^{2(1-s)}\|n_{+1}\|_{L^\infty_tL^2_x}\|u_2\|_{L^2_tL^\infty_x}\|n_{+3}\|_{L^\infty_tL^2_x}\|u_4\|_{L^2_tL^\infty_x}\nonumber\\
&\lesssim&(\frac{N_3}{N})^{2(1-s)}\|n_{+1}\|_{X^{0,\frac{1}{2}+}_+}\|u_2\|_{X^{0,\frac{1}{2}+}}\|n_{+3}\|_{X^{0,\frac{1}{2}+}_+}\|u_4\|_{X^{0,\frac{1}{2}+}}\nonumber\\
&\lesssim&
(\frac{N_3}{N})^{2(1-s)}\frac{1}{<N_1>^{1-s}<N_2>N_3^{1-s}N_4}\|n_{+1}\|_{X^{1-s,\frac{1}{2}+}_+}\|u_2\|_{X^{1,\frac{1}{2}+}}\|n_{+3}\|_{X^{1-s+,\frac{1}{2}+}_+}\|u_4\|_{X^{1,\frac{1}{2}+}}\nonumber\\
&\lesssim&N_{max}^{-\epsilon}
N^{-2+s+}\|In_+\|_{X^{1-s,\frac{1}{2}+}_+}^2\|Iu\|_{X^{1,\frac{1}{2}+}}^2,\label{4.25}
\eea by H\"{o}lder inequality, (\ref{2.8}) and (\ref{2.10}).


{\bf b} $N_1<<N$, $N_2\gtrsim N$.

Then, $|\frac{m(\xi_1+\xi_2)}{m_1m_2}|\sim|\frac{m_2}{m_2}|\sim1$,
and
\bea II'&\lesssim&\|n_{+1}u_2\|_{L^2_{t,x}}\|n_{+3}u_4\|_{L^2_{t,x}}\nonumber\\
&\lesssim&\frac{1}{<N_1>^{1-s}N_2N_3^{1-s}N_4}\|n_{+1}\|_{X^{1-s,\frac{1}{2}+}_+}\|u_2\|_{X^{1,\frac{1}{2}+}}\|n_{+3}\|_{X^{1-s+,\frac{1}{2}+}_+}\|u_4\|_{X^{1,\frac{1}{2}+}}\nonumber\\
&\lesssim&
N_{max}^{-\epsilon}N^{-4+s+}\|In_+\|_{X^{1-s,\frac{1}{2}+}_+}^2\|Iu\|_{X^{1,\frac{1}{2}+}}^2.\label{4.26}
\eea


{\bf c} $N_1\gtrsim N$, $N_2<<N$.

$|\frac{m(\xi_1+\xi_2)}{m_1m_2}|$ is also $\sim 1$.

\bea II'&\lesssim&\|n_{+1}u_2\|_{L^2_{t,x}}\|n_{+3}u_4\|_{L^2_{t,x}}\nonumber\\
&\lesssim&\frac{1}{N_1^{1-s}<N_2>N_3^{1-s}N_4}\|n_{+1}\|_{X^{1-s,\frac{1}{2}+}_+}\|u_2\|_{X^{1,\frac{1}{2}+}}\|n_{+3}\|_{X^{1-s+,\frac{1}{2}+}_+}\|u_4\|_{X^{1,\frac{1}{2}+}}\nonumber\\
&\lesssim&
N_{max}^{-\epsilon}N^{-3+2s+}\|In_+\|_{X^{1-s,\frac{1}{2}+}_+}^2\|Iu\|_{X^{1,\frac{1}{2}+}}^2.\label{4.27}
\eea


{\bf d} $N_1$, $N_2\gtrsim N$.

It exists
$|\frac{m(\xi_1+\xi_2)}{m_1m_2}|\lesssim\frac{1}{m_1m_2}\lesssim(\frac{N_1}{N})^{1-s}(\frac{N_2}{N})^{1-s}$,
thus \bea II'&\lesssim&(\frac{N_1}{N})^{1-s}(\frac{N_2}{N})^{1-s}\|n_{+1}u_2\|_{L^2_{t,x}}\|n_{+3}u_4\|_{L^2_{t,x}}\nonumber\\
&\lesssim&(\frac{N_1}{N})^{1-s}(\frac{N_2}{N})^{1-s}\frac{1}{N_1^{1-s}N_2N_3^{1-s}N_4}\|n_{+1}\|_{X^{1-s,\frac{1}{2}+}_+}\|u_2\|_{X^{1,\frac{1}{2}+}}\|n_{+3}\|_{X^{1-s,\frac{1}{2}+}_+}\|u_4\|_{X^{1,\frac{1}{2}+}}\nonumber\\
&\lesssim&
N_{max}^{-\epsilon}N^{-4+2s+}\|In_+\|_{X^{1-s,\frac{1}{2}+}_+}^2\|Iu\|_{X^{1,\frac{1}{2}+}}^2.\label{4.28}
\eea


{\bf Case 2} $N_{max}\sim N_1\sim N_3\gtrsim N$.

{\bf a} $N_2$, $N_4<<N$.

In this case,
$|\frac{m_{\xi_1+\xi_2}}{m_1m_2}|\lesssim|\frac{m_1}{m_1}|\sim1$,
and
$|\frac{m(\xi_3+\xi_4)-m_3m_4}{m_3m_4}|\lesssim|\frac{m(\xi_3+\xi_4)-m_3}{m_3}|\lesssim|\frac{\nabla
m_3\xi_4}{m_3}|\lesssim\frac{N_4}{N_3}$.

Therefore, \bea II'&\lesssim&\frac{N_4}{N_3}\|n_{+1}u_2\|_{L^2_{t,x}}\|n_{+3}u_4\|_{L^2_{t,x}}\nonumber\\
                   &\lesssim&\frac{N_4}{N_3}\frac{1}{N_1^{1-s}<N_2>N_3^{1-s}<N_4>}\|n_{+1}\|_{X^{1-s,\frac{1}{2}+}_+}\|u_2\|_{X^{1,\frac{1}{2}+}}\|n_{+3}\|_{X^{1-s,\frac{1}{2}+}_+}\|u_4\|_{X^{1,\frac{1}{2}+}}\nonumber\\
                   &\lesssim&N_{max}^{-\epsilon}N^{-3+2s+}\|In_+\|_{X^{1-s,\frac{1}{2}+}_+}^2\|Iu\|_{X^{1,\frac{1}{2}+}}^2.\label{4.29}
\eea


{\bf b} $N_2<<N$, $N_4\gtrsim N$.

So $|\frac{m(\xi_1+\xi_2)}{m_1m_2}|\sim1$,
$|\frac{m(\xi_3+\xi_4)-m_3m_4}{m_3m_4}|\lesssim|\frac{m(\xi_3+\xi_4)}{m_3m_4}|\sim|\frac{m(\xi_1+\xi_2)}{m_3m_4}|\lesssim|\frac{m_1}{m_3m_4}|\sim|\frac{1}{m_4}|\lesssim(\frac{N_4}{N})^{1-s}$,
and \bea II'&\lesssim&(\frac{N_4}{N})^{1-s}\|n_{+1}u_2\|_{L^2_{t,x}}\|n_{+3}u_4\|_{L^2_{t,x}}\nonumber\\
                   &\lesssim&(\frac{N_4}{N})^{1-s}\frac{1}{N_1^{1-s}<N_2>N_3^{1-s}N_4}\|n_{+1}\|_{X^{1-s,\frac{1}{2}+}_+}\|u_2\|_{X^{1,\frac{1}{2}+}}\|n_{+3}\|_{X^{1-s,\frac{1}{2}+}_+}\|u_4\|_{X^{1,\frac{1}{2}+}}\nonumber\\
                   &\lesssim&N_{max}^{-\epsilon}N^{-3+2s+}\|In_+\|_{X^{1-s,\frac{1}{2}+}_+}^2\|Iu\|_{X^{1,\frac{1}{2}+}}^2.\label{4.30}
                   \eea


{\bf c} $N_2\gtrsim N$, $N_4<<N$.

Then
$|\frac{m(\xi_1+\xi_2)}{m_1m_2}|\sim|\frac{m(\xi_3+\xi_4)}{m_1m_2}|\lesssim|\frac{m_3}{m_1m_2}|\sim|\frac{1}{m_2}|\lesssim(\frac{N_2}{N})^{1-s}$,
$|\frac{m(\xi_3+\xi_4)-m_3m_4}{m_3m_4}|\lesssim\frac{N_4}{N_3}$, and
\bea II'&\lesssim&(\frac{N_2}{N})^{1-s}(\frac{N_4}{N_3})\|n_{+1}u_2\|_{L^2_{t,x}}\|n_{+3}u_4\|_{L^2_{t,x}}\nonumber\\
                   &\lesssim&(\frac{N_2}{N})^{1-s}\frac{N_4}{N_3}\frac{1}{N_1^{1-s}N_2N_3^{1-s}<N_4>}\|n_{+1}\|_{X^{1-s,\frac{1}{2}+}_+}\|u_2\|_{X^{1,\frac{1}{2}+}}\|n_{+3}\|_{X^{1-s,\frac{1}{2}+}_+}\|u_4\|_{X^{1,\frac{1}{2}+}}\nonumber\\
                   &\lesssim&N_{max}^{-\epsilon}N^{-4+2s+}\|In_+\|_{X^{1-s,\frac{1}{2}+}_+}^2\|Iu\|_{X^{1,\frac{1}{2}+}}^2.\label{4.31}
\eea


{\bf d} $N_2$, $N_4\gtrsim N$.

{\bf d1} At least one of $N_2$ and $N_4$ $\sim N_1\sim N_3$, w.l.o.g
we suppose $N_2\sim N_1\sim N_3$.

Hence,
$|\frac{m(\xi_1+\xi_2)}{m_1m_2}|\lesssim\frac{1}{m_1m_2}\lesssim(\frac{N_1}{N})^{1-s}(\frac{N_2}{N})^{1-s}$
and
$|\frac{m(\xi_3+\xi_4)-m_3m_4}{m_3m_4}|\lesssim\frac{m(\xi_3+\xi_4)}{m_3m_4}\lesssim\frac{1}{m_3m_4}\lesssim(\frac{N_3}{N})^{1-s}(\frac{N_4}{N})^{1-s}$,
and \begin{eqnarray*}
II'&\lesssim&(\frac{N_1}{N})^{1-s}(\frac{N_2}{N})^{1-s}(\frac{N_3}{N})^{1-s}(\frac{N_4}{N})^{1-s}\frac{1}{N_1^{1-s}}\frac{1}{N_2}\frac{1}{N_3^{1-s}}\frac{1}{N_4}\|In_+\|_{X^{1-s,\frac{1}{2}+}_+}^2\|Iu\|_{X^{1,\frac{1}{2}+}}^2\\
   &\lesssim&N_{max}^{-\epsilon}N^{-4+2s+}\|In_+\|_{X^{1-s,\frac{1}{2}+}_+}^2\|Iu\|_{X^{1,\frac{1}{2}+}}^2.
\end{eqnarray*}

{\bf d2} $N_2$, $N_4<<N_1\sim N_3$.

Then
$\frac{m(\xi_1+\xi_2)}{m_1m_2}\lesssim\frac{m_1}{m_1m_2}\sim\frac{1}{m_2}\lesssim(\frac{N_2}{N})^{1-s}$
and
$|\frac{m(\xi_3+\xi_4)-m_3m_4}{m_3m_4}|\lesssim\frac{m(\xi_3+\xi_4)}{m_3m_4}\lesssim\frac{m_3}{m_3m_4}\sim\frac{1}{m_4}\lesssim(\frac{N_4}{N})^{1-s}$.

\begin{eqnarray*}
II'&\lesssim&(\frac{N_2}{N})^{1-s}(\frac{N_4}{N})^{1-s}\frac{1}{N_1^{1-s}}\frac{1}{N_2}\frac{1}{N_3^{1-s}}\frac{1}{N_4}\|In_+\|_{X^{1-s,\frac{1}{2}+}_+}^2\|Iu\|_{X^{1,\frac{1}{2}+}}^2\\
   &\lesssim&N_{max}^{-\epsilon}N^{-4+2s+}\|In_+\|_{X^{1-s,\frac{1}{2}+}_+}^2\|Iu\|_{X^{1,\frac{1}{2}+}}^2.
\end{eqnarray*}


{\bf Case 3} $N_{max}\sim N_2\sim N_3\gtrsim N$.

{\bf a} $N_1$, $N_4<<N$.

$|\frac{m(\xi_1+\xi_2)}{m_1m_2}|\sim\frac{m_2}{m_2}\sim1$,
$|\frac{m(\xi_3+\xi_4)-m_3m_4}{m_3m_4}|\lesssim\frac{N_4}{N_3}$, and
\bea II'&\lesssim&\frac{N_4}{N_3}\|n_{+1}u_2\|_{L^2_{t,x}}\|n_{+3}u_4\|_{L^2_{t,x}}\nonumber\\
   &\lesssim&\frac{N_4}{N_3}\frac{1}{<N_1>^{1-s}N_2N_3^{1-s}<N_4>}\|n_{+1}\|_{X^{1-s,\frac{1}{2}+}_+}\|u_2\|_{X^{1,\frac{1}{2}+}}\|n_{+3}\|_{X^{1-s,\frac{1}{2}+}_+}\|u_4\|_{X^{1,\frac{1}{2}+}}\nonumber\\
   &\lesssim&N_{max}^{-\epsilon}N^{-3+s+}\|In_+\|_{X^{1-s,\frac{1}{2}+}_+}^2\|Iu\|_{X^{1,\frac{1}{2}+}}^2.\label{4.33}
\eea


{\bf b} $N_1<<N$, $N_4\gtrsim N$.
$|\frac{m(\xi_1+\xi_2)}{m_1m_2}|\sim\frac{m_2}{m_2}\sim1$,
$|\frac{m(\xi_3+\xi_4)-m_3m_4}{m_3m_4}|\lesssim|\frac{m(\xi_3+\xi_4)}{m_3m_4}|\sim|\frac{m(\xi_1+\xi_2)}{m_3m_4}|\lesssim\frac{m(\xi_1+\xi_2)}{m_3m_4}\lesssim\frac{m_2}{m_3m_4}\sim\frac{1}{m_4}\lesssim(\frac{N_4}{N})^{1-s}$
and \bea II'&\lesssim&(\frac{N_4}{N})^{1-s}\|n_{+1}u_2\|_{L^2_{t,x}}\|n_{+3}u_4\|_{L^2_{t,x}}\nonumber\\
   &\lesssim&(\frac{N_4}{N})^{1-s}\frac{1}{<N_1>^{1-s}N_2N_3^{1-s}N_4}\|n_{+1}\|_{X^{1-s,\frac{1}{2}+}_+}\|u_2\|_{X^{1,\frac{1}{2}+}}\|n_{+3}\|_{X^{1-s,\frac{1}{2}+}_+}\|u_4\|_{X^{1,\frac{1}{2}+}}\nonumber\\
   &\lesssim&N_{max}^{-\epsilon}N^{-3+s+}\|In_+\|_{X^{1-s,\frac{1}{2}+}_+}^2\|Iu\|_{X^{1,\frac{1}{2}+}}^2.\label{4.34}
   \eea


{\bf c} $N_1\gtrsim N$, $N_4<<N$.

$|\frac{m(\xi_1+\xi_2)}{m_1m_2}|\sim|\frac{m(\xi_3+\xi_4)}{m_1m_2}|\lesssim\frac{m_3}{m_1m_2}\lesssim\frac{1}{m_1}\lesssim(\frac{N_1}{N})^{1-s}$
and $|\frac{m(\xi_3+\xi_4)-m_3m_4}{m_3m_4}|\lesssim|\frac{\nabla
m_3|\xi_4|}{m_3}|\lesssim\frac{N_4}{N_3}$.

Then \bea II'&\lesssim&(\frac{N_1}{N})^{1-s}\frac{N_4}{N_3}\|n_{+1}u_2\|_{L^2_{t,x}}\|n_{+3}u_4\|_{L^2_{t,x}}\nonumber\\
   &\lesssim&(\frac{N_1}{N})^{1-s}\frac{N_4}{N_3}\frac{1}{N_1^{1-s}N_2N_3^{1-s}<N_4>}\|n_{+1}\|_{X^{1-s,\frac{1}{2}+}_+}\|u_2\|_{X^{1,\frac{1}{2}+}}\|n_{+3}\|_{X^{1-s,\frac{1}{2}+}_+}\|u_4\|_{X^{1,\frac{1}{2}+}}\nonumber\\
   &\lesssim&N_{max}^{-\epsilon}N^{-4+2s+}\|In_+\|_{X^{1-s,\frac{1}{2}+}_+}^2\|Iu\|_{X^{1,\frac{1}{2}+}}^2.\label{4.35}
   \eea


{\bf d} $N_1$, $N_4\gtrsim N$.

$|\frac{m(\xi_1+\xi_2)}{m_1m_2}|\lesssim(\frac{N_1}{N})^{1-s}(\frac{N_2}{N})^{1-s}$,
$|\frac{m(\xi_3+\xi_4)-m_3m_4}{m_3m_4}|\lesssim(\frac{N_3}{N})^{1-s}(\frac{N_4}{N})^{1-s}$,
and \bea II'&\lesssim&(\frac{N_1}{N})^{1-s}(\frac{N_2}{N})^{1-s}(\frac{N_3}{N})^{1-s}(\frac{N_4}{N})^{1-s}\|n_{+1}u_2\|_{L^2_{t,x}}\|n_{+3}u_4\|_{L^2_{t,x}}\nonumber\\
   &\lesssim&(\frac{N_1}{N})^{1-s}(\frac{N_2}{N})^{1-s}(\frac{N_3}{N})^{1-s}(\frac{N_4}{N})^{1-s}\frac{1}{N_1^{1-s}N_2N_3^{1-s}N_4}\|n_{+1}\|_{X^{1-s,\frac{1}{2}+}_+}\|u_2\|_{X^{1,\frac{1}{2}+}}\|n_{+3}\|_{X^{1-s,\frac{1}{2}+}_+}\|u_4\|_{X^{1,\frac{1}{2}+}}\nonumber\\
   &\lesssim&N_{max}^{-\epsilon}N^{-4+2s+}\|In_+\|_{X^{1-s,\frac{1}{2}+}_+}^2\|Iu\|_{X^{1,\frac{1}{2}+}}^2.\label{4.36}
\eea


{\bf Case 4} $N_{max}\sim N_1\sim N_4\gtrsim N$.

{\bf a} $N_2$, $N_3<<N$.

$|\frac{m(\xi_1+\xi_2)}{m_1m_2}|\lesssim\frac{m_1}{m_1}\sim1$, and
$|\frac{m(\xi_3+\xi_4)-m_3m_4}{m_3m_4}|\lesssim|\frac{\nabla
m_4\xi_3}{m_4}|\lesssim\frac{N_3}{N_4}$.

Thus \bea II'&\lesssim&\frac{N_3}{N_4}\|n_{+1}u_2\|_{L^2_{t,x}}\|n_{+3}u_4\|_{L^2_{t,x}}\nonumber\\
   &\lesssim&\frac{N_3}{N_4}\frac{1}{N_1^{1-s}<N_2><N_3>^{1-s}N_4}\|n_{+1}\|_{X^{1-s,\frac{1}{2}+}_+}\|u_2\|_{X^{1,\frac{1}{2}+}}\|n_{+3}\|_{X^{1-s,\frac{1}{2}+}_+}\|u_4\|_{X^{1,\frac{1}{2}+}}\nonumber\\
   &\lesssim&N_{max}^{-\frac{\epsilon}{2}}N^{-3+2s+}\|In_+\|_{X^{1-s,\frac{1}{2}+}_+}^2\|Iu\|_{X^{1,\frac{1}{2}+}}^2.\label{4.37}
\eea


{\bf b} $N_2<<N$, $N_3\gtrsim N$.

$|\frac{m(\xi_1+\xi_2)}{m_1m_2}|\sim1$,
$|\frac{m(\xi_3+\xi_4)-m_3m_4}{m_3m_4}|\lesssim\frac{m(\xi_3+\xi_4)}{m_3m_4}\sim\frac{m(\xi_1+\xi_2)}{m_3m_4}\lesssim\frac{m_1}{m_3m_4}\sim\frac{1}{m_3}\lesssim(\frac{N_3}{N})^{1-s}$,
and \bea II'&\lesssim&(\frac{N_3}{N})^{1-s}\|n_{+1}u_2\|_{L^2_{t,x}}\|n_{+3}u_4\|_{L^2_{t,x}}\nonumber\\
   &\lesssim&(\frac{N_3}{N})^{1-s}\frac{1}{N_1^{1-s}<N_2>N_3^{1-s}N_4}\|n_{+1}\|_{X^{1-s,\frac{1}{2}+}_+}\|u_2\|_{X^{1,\frac{1}{2}+}}\|n_{+3}\|_{X^{1-s,\frac{1}{2}+}_+}\|u_4\|_{X^{1,\frac{1}{2}+}}\nonumber\\
   &\lesssim&N_{max}^{-\epsilon}N^{-3+2s+}\|In_+\|_{X^{1-s,\frac{1}{2}+}_+}^2\|Iu\|_{X^{1,\frac{1}{2}+}}^2.\label{4.38}
\eea


{\bf c} $N_2\gtrsim N$, $N_3<<N$.

$|\frac{m(\xi_1+\xi_2)}{m_1m_2}|\sim\frac{m(\xi_3+\xi_4)}{m_1m_2}\lesssim\frac{m_4}{m_1m_2}\sim\frac{1}{m_2}\lesssim(\frac{N_2}{N})^{1-s}$,
$|\frac{m(\xi_3+\xi-4)-m_3m_4}{m_3m_4}|\lesssim\frac{N_3}{N_4}$.

Then \bea II'&\lesssim&(\frac{N_2}{N})^{1-s}\frac{N_3}{N_4}\|n_{+1}u_2\|_{L^2_{t,x}}\|n_{+3}u_4\|_{L^2_{t,x}}\nonumber\\
   &\lesssim&(\frac{N_2}{N})^{1-s}\frac{N_3}{N_4}\frac{1}{N_1^{1-s}N_2<N_3>^{1-s}N_4}\|n_{+1}\|_{X^{1-s,\frac{1}{2}+}_+}\|u_2\|_{X^{1,\frac{1}{2}+}}\|n_{+3}\|_{X^{1-s,\frac{1}{2}+}_+}\|u_4\|_{X^{1,\frac{1}{2}+}}\nonumber\\
   &\lesssim&N_{max}^{-\frac{\epsilon}{2}}N^{-4+2s+}\|In_+\|_{X^{1-s,\frac{1}{2}+}_+}^2\|Iu\|_{X^{1,\frac{1}{2}+}}^2.\label{4.39}
\eea


{\bf d} $N_2$, $N_3\gtrsim N$.

The same as Case 3(d).


{\bf Case 5} $N_{max}\sim N_2\sim N_4\gtrsim N$.

{\bf a} $N_1$, $N_3<<N$.

$|\frac{m(\xi_1+\xi_2)}{m_1m_2}|\lesssim1$,
$|\frac{m(\xi_3+\xi_4)-m_3m_4}{m_3m_4}|\lesssim\frac{N_3}{N_4}$, and
\bea II'&\lesssim&\frac{N_3}{N_4}\|n_{+1}u_2\|_{L^2_{t,x}}\|n_{+3}u_4\|_{L^2_{t,x}}\nonumber\\
   &\lesssim&\frac{N_3}{N_4}\frac{1}{<N_1>^{1-s}N_2<N_3>^{1-s}N_4}\|n_{+1}\|_{X^{1-s,\frac{1}{2}+}_+}\|u_2\|_{X^{1,\frac{1}{2}+}}\|n_{+3}\|_{X^{1-s,\frac{1}{2}+}_+}\|u_4\|_{X^{1,\frac{1}{2}+}}\nonumber\\
   &\lesssim&N_{max}^{-\frac{\epsilon}{2}}N^{-3+s+}\|In_+\|_{X^{1-s,\frac{1}{2}+}_+}^2\|Iu\|_{X^{1,\frac{1}{2}+}}^2.\label{4.40}
\eea


{\bf b} $N_1<<N$, $N_3\gtrsim N$.

$|\frac{m(\xi_1+\xi_2)}{m_1m_2}|\lesssim1$,
$|\frac{m(\xi_3+\xi_4)-m_3m_4}{m_3m_4}|\lesssim\frac{m(\xi_3+\xi_4)}{m_3m_4}\sim\frac{m(\xi_1+\xi_2)}{m_3m_4}\lesssim\frac{m_2}{m_3m_4}\sim\frac{1}{m_3}\lesssim(\frac{N_3}{N})^{1-s}$,
and \bea II'&\lesssim&(\frac{N_3}{N})^{1-s}\|n_{+1}u_2\|_{L^2_{t,x}}\|n_{+3}u_4\|_{L^2_{t,x}}\nonumber\\
   &\lesssim&(\frac{N_3}{N})^{1-s}\frac{1}{<N_1>^{1-s}N_2N_3^{1-s}N_4}\|n_{+1}\|_{X^{1-s,\frac{1}{2}+}_+}\|u_2\|_{X^{1,\frac{1}{2}+}}\|n_{+3}\|_{X^{1-s,\frac{1}{2}+}_+}\|u_4\|_{X^{1,\frac{1}{2}+}}\nonumber\\
   &\lesssim&N_{max}^{-\epsilon}N^{-3+s+}\|In_+\|_{X^{1-s,\frac{1}{2}+}_+}^2\|Iu\|_{X^{1,\frac{1}{2}+}}^2.\label{4.41}
\eea


{\bf c} $N_1\gtrsim N$, $N_3<<N$.

$|\frac{m(\xi_1+\xi_2)}{m_1m_2}|\sim\frac{m(\xi_3+\xi_4)}{m_1m_2}\lesssim\frac{m_4}{m_1m_2}\sim\frac{1}{m_1}\lesssim(\frac{N_1}{N})^{1-s}$,
$|\frac{m(\xi_3+\xi_4)-m_3m_4}{m_3m_4}|\lesssim\frac{N_3}{N_4}$, and

\bea II'&\lesssim&(\frac{N_1}{N})^{1-s}\frac{N_3}{N_4}\|n_{+1}u_2\|_{L^2_{t,x}}\|n_{+3}u_4\|_{L^2_{t,x}}\nonumber\\
   &\lesssim&(\frac{N_1}{N})^{1-s}\frac{N_3}{N_4}\frac{1}{N_1^{1-s}N_2<N_3>^{1-s}N_4}\|n_{+1}\|_{X^{1-s,\frac{1}{2}+}_+}\|u_2\|_{X^{1,\frac{1}{2}+}}\|n_{+3}\|_{X^{1-s,\frac{1}{2}+}_+}\|u_4\|_{X^{1,\frac{1}{2}+}}\nonumber\\
   &\lesssim&N_{max}^{-\frac{\epsilon}{2}}N^{-4+2s+}\|In_+\|_{X^{1-s,\frac{1}{2}+}_+}^2\|Iu\|_{X^{1,\frac{1}{2}+}}^2.\label{4.42}
\eea


{\bf d} $N_1$, $N_3\gtrsim N$.

The same as Case 2(d).

\end{proof}


Finally, let's consider $III$.

\begin{lem}\label{lem3}
\be\label{4.43} III'\lesssim N^{-2+s}\delta^{\frac{1}{2}-}
\|In_+\|_{X^{1-s+,\frac{1}{2}+}_+}\|Iu\|_{X^{1,\frac{1}{2}+}}^2. \ee
\end{lem}
\begin{proof}

With the same notations and argument as before, to prove this lemma,
we just need the following estimate:
\be\label{4.44}III'=N_1\int_0^\delta\int_*|\frac{m(\xi_2+\xi_3)-m_2m_3}{m_2m_3}|\hat{n}_{+1}\hat{u}_2\hat{u}_3\lesssim
N^{-2+s+}\delta^{\frac{1}{2}-}\|In_+\|_{X^{1-s+,\frac{1}{2}+}_+}\|Iu\|_{X^{1,\frac{1}{2}+}}^2.\ee

If both $N_2$ and $N_3<<N$, l.h.s would be zero, and the inequality
holds.

On the other hand, w.l.o.g we assume $N_3\leqslant N_2$.

So let's suppose $N_2\gtrsim N$. Since $\sum_{i=1}^3\xi_i=0$,
$N_1\lesssim N_2.$

Now, we'll discuss in two subcases.

{\bf Case 1} $N_2\gtrsim N>>N_3$.

As $\sum_{i=1}^3\xi_i=0$, then $N_1\sim N_2$, and
$|\frac{m(\xi_2+\xi_3)-m_2m_3}{m_2m_3}|\sim|\frac{m(\xi_2+\xi_3)-m_2}{m_2}|\lesssim|\frac{\nabla
m_2\xi_3}{m_2}|\lesssim\frac{N_3}{N_2}$.

\bea
III'&\lesssim&N_1\frac{N_3}{N_2}\|n_{+1}\|_{L^2_{t,x}}\|u_2u_3\|_{L^2_{t,x}}\nonumber\\
    &\lesssim&N_1\frac{N_3}{N_2}\|n_{+1}\|_{X^{0,0}_+}(\frac{N_3}{N_2})^{\frac{1}{2}}\|u_2\|_{X^{0,\frac{1}{2}+}}\|u_3\|_{X^{0,\frac{1}{2}+}}\nonumber\\
    &\lesssim&N_1\frac{N_3}{N_2}(\frac{N_3}{N_2})^\frac{1}{2}\frac{1}{N_1^{1-s}}\frac{1}{N_2}\frac{1}{<N_3>}\delta^{\frac{1}{2}-}\|n_{+1}\|_{X^{1-s,0}_+}\|u_2\|_{X^{1,\frac{1}{2}+}}\|u_3\|_{X^{1,\frac{1}{2}+}}\nonumber\\
    &\lesssim&N_{max}^{-\frac{\epsilon}{2}}N^{-2+s+}\delta^{\frac{1}{2}-}\|In_+\|_{X^{1-s,\frac{1}{2}+}_+}\|Iu\|_{X^{1,\frac{1}{2}+}}^2,\label{4.45}
\eea by H\"{o}lder inequality, (\ref{2.7}) and (\ref{2.11}).


{\bf Case 2} $N_2\geqslant N_3\gtrsim N$.

{\bf Subcase a} $N_1\sim N_2\geqslant N_3\gtrsim N$.

$|\frac{m(\xi_2+\xi_3)-m_2m_3}{m_2m_3}|\lesssim\frac{m(\xi_2+\xi_3)}{m_2m_3}\sim\frac{m_1}{m_2m_3}\sim\frac{1}{m_3}\lesssim(\frac{N_3}{N})^{1-s}$,
then \bea
III'&\lesssim&N_1(\frac{N_3}{N})^{1-s}\|n_{+1}\|_{L^2_{t,x}}\|u_2u_3\|_{L^2_{t,x}}\nonumber\\
    &\lesssim&N_1(\frac{N_3}{N})^{1-s}\|n_{+1}\|_{X^{0,0}_+}(\frac{N_3}{N_2})^{\frac{1}{2}}\|u_2\|_{X^{0,\frac{1}{2}+}}\|u_3\|_{X^{0,\frac{1}{2}+}}\nonumber\\
    &\lesssim&N_1(\frac{N_3}{N})^{1-s}(\frac{N_3}{N_2})^\frac{1}{2}\frac{1}{N_1^{1-s}}\frac{1}{N_2}\frac{1}{N_3}\delta^{\frac{1}{2}-}\|n_{+1}\|_{X^{1-s,0}_+}\|u_2\|_{X^{1,\frac{1}{2}+}}\|u_3\|_{X^{1,\frac{1}{2}+}}\nonumber\\
    &\lesssim&N_{max}^{-\epsilon}N^{-2+s+}\delta^{\frac{1}{2}-}\|In_+\|_{X^{1-s,\frac{1}{2}+}_+}\|Iu\|_{X^{1,\frac{1}{2}+}}^2,\label{4.46}
\eea with the same reason as case 1.


{\bf Subcase b} $N_2\geqslant N_3\gtrsim N$, $N_2>>N_1$.

So $N_2\sim N_3$ and
$|\frac{m(\xi_2+\xi_3)-m_2m_3}{m_2m_3}|\lesssim\frac{m(\xi_2+\xi_3)}{m_2m_3}\lesssim\frac{1}{m_2m_3}\sim\frac{1}{m_2^2}\lesssim(\frac{N_2}{N})^{2(1-s)}$.
\bea
III'&\lesssim&N_1(\frac{N_2}{N})^{2(1-s)}\|n_{+1}\|_{L^2_{t,x}}\|u_2u_3\|_{L^2_{t,x}}\nonumber\\
    &\lesssim&N_1(\frac{N_2}{N})^{2(1-s)}\|_{X^{0,0}_+}(\frac{N_3}{N_2})^{\frac{1}{2}}\|u_2\|_{X^{0,\frac{1}{2}+}}\|u_3\|_{X^{0,\frac{1}{2}+}}\nonumber\\
    &\lesssim&N_1(\frac{N_2}{N})^{2(1-s)}(\frac{N_3}{N_2})^\frac{1}{2}\frac{1}{<N_1>^{1-s}}\frac{1}{N_2}\frac{1}{N_3}\delta^{\frac{1}{2}-}\|n_{+1}\|_{X^{1-s,0}_+}\|u_2\|_{X^{1,\frac{1}{2}+}}\|u_3\|_{X^{1,\frac{1}{2}+}}\nonumber\\
    &\lesssim&N_{max}^{-\epsilon}N^{-2+s+2\epsilon}\delta^{\frac{1}{2}-}\|In_+\|_{X^{1-s,\frac{1}{2}+}_+}\|Iu\|_{X^{1,\frac{1}{2}+}}^2.\label{4.47}
\eea
\end{proof}


Now, combing the results of three lemmas above, we can get
Proposition \ref{prop4.1} easily.

\end{proof}


\section{Proof of Theorem \ref{th1}}

Let \be\label{5.1}\Sigma_u(t)=\sup_{0\leqslant\tau\leqslant
t}\|I_N<\Lambda>u(\tau)\|_{L^2},\ee \be\label{5.2}
\Sigma_{n_+}(t)=\sup_{0\leqslant\tau\leqslant
t}\|I_N<\Lambda>^{1-s}n_+(\tau)\|_{L^2}\ee
\be\label{5.3}\widetilde{\Sigma}_u(t)=\sup_{0\leqslant\tau\leqslant
t}\|I_Nu(\tau)\|_{X^{1,\frac{1}{2}+}},\ee\be\label{5.4}\widetilde{\Sigma}_{n_+}(t)=\sup_{0\leqslant\tau\leqslant
t}\|I_Nn_+(\tau)\|_{X^{1-s,\frac{1}{2}+}_+}\ee and
\be\label{5.5}\Lambda(t)=\sup_{0\leqslant\tau\leqslant
t}\|<\Lambda>^su(\tau)\|_{L^2}\ee

First of all, we'll prove the following proposition.
\begin{prop}\label{prop5.1} With the condition of Theorem \ref{th1}
for $1>s>\frac{16}{17}$, $\forall T<T^*<\infty$ and close to $T^*$
enough, \be\label{5.6}
|H(T)|=|H(Iu(T),In_+(T))|\lesssim\Lambda^{p(s)},\ee where
$N\sim\Lambda(T)^{\frac{10+34\epsilon}{7s-6-(35-34s)\epsilon}}$,
$\epsilon$ small enough such that $0<\epsilon<\frac{17-16s}{69-68s}$
and $p(s)<2$.
\end{prop}
\begin{proof}
From Proposition \ref{prop4.1} and the condition of Theorem
\ref{th1}, there exists \bea |H(\delta)-H(0)|&\lesssim&
N^{-2+s+}\delta^{0+}\|In_+\|_{X^{1-s,\frac{1}{2}+}_+}\|Iu\|_{X^{1,\frac{1}{2}+}}^2+N^{-2+s+}\|In_+\|_{X^{1-s,\frac{1}{2}+}_+}^2\|Iu\|_{X^{1,\frac{1}{2}+}}^2\nonumber\\
&\lesssim&N^{-2+s+}\|In_+\|_{X^{1-s,\frac{1}{2}+}_+}^2\|Iu\|_{X^{1,\frac{1}{2}+}}^2,\label{5.7}\eea
since
$\|In_+\|_{X^{1-s,\frac{1}{2}+}_+}+\|Iu\|_{X^{1,\frac{1}{2}+}}\gtrsim
\|In_+\|_{H^{1-s}}+\|Iu\|_{H^1}\rightarrow\infty$ for $t\rightarrow
T^*$.

On the other hand, by Proposition \ref{prop3.1}, we choose
$\delta^{-1}\sim \Sigma_u(T)^{2+17\epsilon}+(\sup_{0\leqslant
\tau\leqslant
T}\frac{\|Iu(\tau)\|_{H^1}^2}{\|In_+(\tau)\|_{H^{1-s}}})^{2+17\epsilon}$.
As from (\ref{3.14}), it has $\|In_+(t)\|_{H^{1-s}}\geqslant c>0$,
as $t\rightarrow T^*$. w.l.o.g we suppose $\|In_+(t)\|\geqslant c$
for $0\leqslant t<T^*$, otherwise, we just need to calculate from
$H(t^*)$ for some $t^*<T^*$. Thus $\delta^{-1}\lesssim
\Sigma_{n_+}(T)^{2+17\epsilon}+\Sigma_u(T)^{2(2+17\epsilon)}$, and
the number of iteration steps to reach the given time $T$ is
$\frac{T}{\delta}\lesssim
T(\Sigma_{n_+}(T)^{2+17\epsilon}+\Sigma_u(T)^{2(2+17\epsilon)})$.

Combining these estimates with (\ref{5.7}), the whole increment of
energy is \bea
&&T(\Sigma_{n_+}(T)^{2+17\epsilon}+\Sigma_u(T)^{2(2+17\epsilon)})N^{-2+s+\epsilon}\widetilde{\Sigma}_{n_+}(T)^2\widetilde{\Sigma}_u(T)^2\nonumber\\
&\lesssim&N^{-2+s+\epsilon}(\widetilde{\Sigma}_{n_+}(T)^{4+17\epsilon}\widetilde{\Sigma}_u(T)^2+\widetilde{\Sigma}_{n_+}(T)^{2}\widetilde{\Sigma}_u(T)^{6+34\epsilon})\label{5.8}.
\eea

Then, from (\ref{3.15}) for $T<T^*$, \be\label{5.9}
\|(In_+)\|_{X^{1-s,\frac{1}{2}+}_+}\leqslant
c\|In_{+0}\|_{H^{1-s}}+cT^{\frac{1}{2}-4\epsilon}\|Iu\|_{X^{1,\frac{1}{2}+}}^2.\ee
Hence, \be\label{5.10} \widetilde{\Sigma}_{n_+}(T)\lesssim
N^{1-s}\|n_{+0}\|_{L^2}+\widetilde{\Sigma}_u(T)^2\lesssim
N^{1-s}+\widetilde{\Sigma}_u(T)^2,\ee then put it into (\ref{5.8}),
and by the relationship (\ref{3.9}), then \bea (\ref{5.8})&\lesssim&
N^{-2+s+\epsilon}(N^{(4+17\epsilon)(1-s)}\widetilde{\Sigma}_u(T)^2+\widetilde{\Sigma}_u(T)^{10+34\epsilon}+N^{2(1-s)}\widetilde{\Sigma}_u(T)^{6+34\epsilon}+\widetilde{\Sigma}_u(T)^{10+34\epsilon})\nonumber\\
&\lesssim&N^{-2+s+\epsilon}(N^{(4+17\epsilon)(1-s)}\Sigma_u(T)^2+N^{2(1-s)}\Sigma_u(T)^{6+34\epsilon}+\Sigma_u(T)^{10+34\epsilon})\nonumber\\
&\lesssim&N^{-2+s+\epsilon}(N^{(4+17\epsilon)(1-s)}N^{2(1-s)}\Lambda(T)^2+N^{2(1-s)}N^{(6+34\epsilon)(1-s)}\Lambda(T)^{6+34\epsilon}+N^{10+34\epsilon(1-s)}\Lambda(T)^{10+34\epsilon})\nonumber\\
&\lesssim&N^{4-5s+(18-17s)\epsilon}\Lambda(T)^2+N^{6-7s+(35-34s)\epsilon}\Lambda(T)^{6+34\epsilon}+N^{8-9s+(35-34s)\epsilon}\Lambda(T)^{10+34\epsilon}\label{5.11}.\eea

On the other hand, \bea |H(0)|=|H(Iu_0,In_{+0})|&=&\|\nabla
Iu_0\|_{L^2}^2+\frac{1}{2}\|In_{+0}\|_{L^2}^2+\frac{1}{2}\int(In_{+0}+\overline{In}_{+0})|Iu|^2dx\nonumber\\
                       &\lesssim&\|\nabla Iu_0\|_{L^2}^2+\|In_{+0}\|_{L^2}^2+\|In_{+0}\|_{L^2}\|Iu_0\|_{L^4}^2\nonumber\\
                       &\lesssim&\|\nabla
                       Iu_0\|_{L^2}^2+\|In_{+0}\|_{L^2}^2+\|Iu_0\|_{L^4}^4\nonumber\\
                       &\lesssim&\|\nabla
                       Iu_0\|_{L^2}^2+\|In_{+0}\|_{L^2}^2+\|Iu_0\|_{L^2}^2\|\nabla Iu_0\|_{L^2}^2\nonumber\\
                       &\lesssim&\|\nabla
                       Iu_0\|_{L^2}^2+\|In_{+0}\|_{L^2}^2+\|u_0\|_{L^2}^2\|\nabla Iu_0\|_{L^2}^2\nonumber\\
                       &\lesssim&\|\nabla
                       Iu_0\|_{L^2}^2+\|In_{+0}\|_{L^2}^2+\|\nabla Iu_0\|_{L^2}^2\nonumber\\
                       &\lesssim&\|\nabla
                       Iu_0\|_{L^2}^2+\|In_{+0}\|_{L^2}^2\nonumber\\
                       &\lesssim&N^{2(1-s)}\|u_{0}\|_{H^1}^2+\|n_{+0}\|_{L^2}^2\lesssim
                       N^{2(1-s)},\label{5.12}\eea
Hence, \bea |H(T)|&\leqslant& |H(0)|+|H(T)-T(0)|\nonumber\\
&\lesssim&
N^{2(1-s)}+N^{4-5s+(18-17s)\epsilon}\Lambda(T)^2+N^{6-7s+(35-34s)\epsilon}\Lambda(T)^{6+34\epsilon}\nonumber\\
&&+N^{8-9s+(35-34s)\epsilon}\Lambda(T)^{10+34\epsilon}\label{5.13}.
\eea

Then, choose
$N=\Lambda^{\frac{10+34\epsilon}{7s-6-(35-34s)\epsilon}}$, so that
the first and fourth terms in (\ref{5.13}) give comparable
contributions. A calculation reveals that the second and third terms
in (\ref{5.13}) produces a smaller correction. Thus \be\label{5.14}
p(s)=2(1-s)\frac{10+34\epsilon}{7s-6-(35-34s)\epsilon}<2\Leftrightarrow
s>\frac{16}{17}\ and\ 0<\epsilon<\frac{17s-16}{69-68s}. \ee
\end{proof}

Now we turn to prove Theorem \ref{th1}.

\begin{proof}

Let $\{t_k\}_{k=1}^\infty$ be a sequence such that $t_k\uparrow
T^*$ as $k\rightarrow\infty$, and for each $t_k$,
\[\|u(t_k)\|_{H^s}=\Lambda(t_k).\] By the result of Corollary
\ref{cor3.2}, that $\|u(t)\|_{H^s}\rightarrow\infty$, it's
achievable.

Denote $u_k=u(t_k)$, and $Iu_k=I_{N(t_k)}u(t_k)$, with $N(t_k)$
taken as in Proposition \ref{prop5.1}.

Then, let $\lambda_k=\|Iu_k\|_{H^1}\geqslant\Lambda(t_k)$. Do the
scaling as follows: \bea
\tilde{u}_k&=&\lambda_k^{-1}Iu(t_k,x\lambda_k^{-1})\label{5.15}\\
\tilde{n}_k&=&\lambda_k^{-2}In(t_k,x\lambda_k^{-1})\label{5.16},\eea
and by direct calculations, we have
\be\label{5.17}\|\tilde{u}_k\|_{L^2}=\|Iu_k\|_{L^2}\leqslant\|u_k\|_{L^2}=\|u_0\|_{L^2},\ee
\be\label{5.18}\|\nabla\tilde{u}_k\|_{L^2}\leqslant1,\ee
\be\label{5.20}\lim_{k\rightarrow\infty}\|\tilde{u}_k\|_{L^2}=1,\ee
    and
\be\label{5.19}\lim_{k\rightarrow\infty}\|\nabla\tilde{u}_k\|_{L^2}=1\ee
since $N(t_k)\rightarrow\infty$ for $t\rightarrow T^*$ by
Proposition \ref{prop5.1}.

Thus, $\{\tilde{u}_k\}_{k=1}^\infty$ is a bounded sequence in $H^1$
and has a weakly convergent subsequence, which we still denote as
$\{\tilde{u}_k\}$, and $\tilde{u}\in H^1$, such that \be\label{5.20}
\tilde{u}_n\rightharpoonup\tilde{u}\ \ \ in\ H^1.\ee Then, as $u$ is
radial, then by Radial Compactness Lemma, it exists
\be\label{5.21}\tilde{u}_k\rightarrow\tilde{u}\ \ \ in\ L^4.\ee

On the other hand, let \be\label{5.22} E(Iu)=\|\nabla
Iu(t)\|_{L^2}^2-\frac{1}{2}\|Iu(t)\|_{L^4}^4,\ee and \be\label{5.23}
H_1(Iu,In)=\|\nabla Iu(t)\|_{L^2}^2+\frac{1}{2}\|In\|_{L^2}^2+ \int
In|Iu|^2=E(Iu)+\frac{1}{2}\int(In+|Iu|^2)^2.\ee Hence
\be\label{5.24}H(t)=H_1(Iu,In)+\frac{1}{2}\|Iv\|_{L^2}^2.\ee

Combing all the above estimates together, we have \bea
E(\tilde{u}_k)&=&\lambda_k^{-2}E(Iu_k)\nonumber\\
                &\leqslant& \lambda_k^{-2}H_1(Iu_k,In_k)\nonumber\\
                &\leqslant& \lambda_k^{-2}H(Iu_k,In_k)\nonumber\\
                &\leqslant& c\lambda_k^{-2}\Lambda^{p(s)}(t_k)\nonumber\\
                &\leqslant& c\Lambda_k^{p(s)-2}\rightarrow 0, \nonumber\eea
as $k\rightarrow\infty$, by Proposition \ref{prop5.1} and the
definition of $t_k$.

Thus,
\be\label{5.25}\limsup_{k\rightarrow\infty}E(\tilde{u}_k)\leqslant0,\ee
and
\be\label{5.26}\limsup_{k\rightarrow\infty}H_1(\tilde{u}_k,\tilde{n}_k)\leqslant0.\ee
Therefore, \be\label{5.27}
\liminf_{k\rightarrow\infty}\|\tilde{u}_k\|_{L^4}^4=2\liminf_{k\rightarrow\infty}\|\nabla\tilde{u}_k\|_{L^2}^2-2E(\tilde{u}_k)\geqslant2,\ee
and
\be\label{5.28}0\geqslant\limsup_{k\rightarrow\infty}H_1(\tilde{u}_k,\tilde{n}_k)\geqslant\frac{1}{2}\limsup_{k\rightarrow\infty}(\int(\tilde{n}_k+|\tilde{u}_k|^2)^2-\|\tilde{u}_k\|_{L^4}^4),\ee
or in other words,
\[0\geqslant\limsup_{k\rightarrow\infty}H_1(\tilde{u}_k,\tilde{n}_k)\geqslant\limsup_{k\rightarrow\infty}(\frac{1}{2}\|\tilde{n}_k\|_{L^2}^2+\int\tilde{n}_k|\tilde{u}_k|^2)\geqslant\limsup_{k\rightarrow\infty}(\frac{1}{2}\|\tilde{n}_k\|_{L^2}^2-\frac{1}{4}\|\tilde{n}\|_{L^2}^2-\|\tilde{u}_k\|_{L^4}^4),\]
i.e.\be\label{5.29}\limsup_{k\rightarrow\infty}\|\tilde{n}_k\|_{L^2}^2\leqslant4\liminf_{k\rightarrow\infty}\|\tilde{u}_k\|_{L^4}^4\leqslant
c,\ee by sobolev embedding theory, (\ref{5.17}) and (\ref{5.18}).

\begin{cl}\label{cl5.1}
 $\forall\ R>0$,
 \be\label{5.30}\liminf_{k\rightarrow\infty}\|Iu_k\|_{L^2(B(0,R))}\geqslant\|Q\|_{L^2}\ee
 and
 \be\label{5.31}\liminf_{k\rightarrow\infty}\|In_k\|_{L^1(B(0,R))}\geqslant
 m_n,\ee where $m_n>0$ depending on the initial data.
\end{cl}
\begin{proof}
If the claim doesn't exist, then there is a subsequence of
$\{t_k\}$, (still denote it as $\{t_k\}$), such that
\be\label{5.32}\limsup_{k\rightarrow\infty}\int_{|x|<R_0}|Iu_k|^2\leqslant\|Q\|_{L^2}^2-\delta_0,\ee
or
\be\label{5.33}\limsup_{k\rightarrow\infty}\int_{|x|<R_0}|In_k|=0,\ee
for some $R_0>0$, and $\delta_0>0$.

Then by scaling, $\forall R>0$,
\be\label{5.34}\limsup_{k\rightarrow\infty}\int_{|x|<R}|\tilde{u}_k|^2\leqslant\|Q\|_{L^2}^2-\delta_0,\ee
or
\be\label{5.35}\limsup_{k\rightarrow\infty}\int_{|x|<R}|\tilde{n}_k|=0,\ee
as $\lambda_k\rightarrow\infty$ for $k\rightarrow\infty$.

From (\ref{5.21}) and (\ref{5.27}), there exists
\be\label{5.36}\|\tilde{u}\|_{L^4}^4\geqslant2.\ee And also by
(\ref{5.20}) and (\ref{5.34}), we have
\[\int_{|x|<R}|\tilde{u}|^2\leqslant\liminf_{k\rightarrow\infty}\int_{|x|<R}|\tilde{u}_k|^2\leqslant\|Q\|_{L^2}^2-\delta_0,\]
for any $R>0$. Hence by letting $R\rightarrow\infty$,
\be\label{5.37}\|\tilde{u}\|_{L^2}^2\leqslant\|Q\|_{L^2}^2-\delta_0.\ee

On the other hand, from (\ref{5.29}), we can see that
$\{\tilde{n}_k\}$ is bounded in $L^2$, hence there is $\tilde{n}$,
such that \be\label{5.38}\tilde{n}_k\rightharpoonup\tilde{n}\ \ \ \
in\ L^2.\ee

From (\ref{5.38}) and (\ref{5.35})
\[\int_{|x|<R}|\tilde{n}|\leqslant cR^{\frac{1}{2}}(\int_{|x|<R}|\tilde{n}|^2)^\frac{1}{2}\leqslant cR^{\frac{1}{2}}\liminf_{k\rightarrow\infty}(\int_{|x|<R}|\tilde{n}_k|^2)^\frac{1}{2}=0.\]
i.e. \be\label{5.39} \tilde{n}=0,\ \  a.e.\ee by letting
$R\rightarrow\infty$.

Therefore,
\be\label{5.40}\|\tilde{u}\|_{L^2}^2\leqslant\|Q\|_{L^2}^2-\delta_0\
\ or\ \ \tilde{n}=0.\ee

Furthermore, since $\tilde{u}_k^2\rightarrow\tilde{u}^2$ and
$\tilde{n}_k\rightharpoonup\tilde{n}$ in $L^2$, we have
\be\label{5.41}\int\tilde{n}_k|\tilde{u}_k|^2\rightarrow\int\tilde{n}|\tilde{u}|^2,\
\ as\ \ k\rightarrow\infty.\ee

Therefore,
\be\label{5.42}H_1(\tilde{u},\tilde{n})=\|\nabla\tilde{u}\|_{L^2}^2+\frac{1}{2}\|\tilde{n}\|_{L^2}^2+\int\tilde{n}|\tilde{u}|^2\leqslant\liminf_{k\rightarrow\infty}(\|\nabla\tilde{u}_k\|_{L^2}^2+\frac{1}{2}\|\tilde{n}_k\|_{L^2}^2+\int\tilde{n}_k|\tilde{u}_k|^2)=\liminf_{k\rightarrow\infty}H_1(\tilde{u}_k,\tilde{n}_k)\leqslant0,
\ee or equivalently,
\be\label{5.43}E(\tilde{u})+\frac{1}{2}\int(\tilde{n}+|\tilde{u}|^2)^2\leqslant0.\ee

{\bf Case 1}\ \ \ \  If
$\|\tilde{u}\|_{L^2}^2\leqslant\|Q\|_{L^2}^2-\delta_0$.

Then, by (\ref{5.43}) and sharp Gagliardo-Nirenberg (\ref{2.14}),
we have \bea 0&\geqslant&
E(\tilde{u})=\|\nabla\tilde{u}\|_{L^2}^2-\frac{1}{2}\|\tilde{u}\|_{L^4}^4\nonumber\\
           &\geqslant&\|\nabla\tilde{u}\|_{L^2}^2-\frac{\|\tilde{u}\|_{L^2}^2}{\|Q\|_{L^2}^2}\|\nabla\tilde{u}\|_{L^2}^2\nonumber\\
           &\geqslant&(1-\frac{\|Q\|_{L^2}^2-\delta_0}{\|Q\|_{L^2}^2})\|\nabla\tilde{u}\|_{L^2}^2\nonumber\\
           &=&\frac{\delta_0}{\|Q\|_{L^2}^2}\|\nabla\tilde{u}\|_{L^2}^2.\label{5.44}
\eea Because of (\ref{5.36}), $\|\nabla\tilde{u}\|_{L^2}^2\neq0$,
which is a contradiction.

{\bf Case 2}\ \ \ \  If $\tilde{n}=0$.

Then \be\label{5.45} 0\geqslant
H_1(\tilde{u},\tilde{n})=\|\nabla\tilde{u}\|_{L^2}^2,\ee which is
also a contradiction.
\end{proof}

With Claim \ref{cl5.1}, we can get the result of the Theorem
quickly.

That is,
\be\label{5.46}\|Q\|_{L^2}\leqslant\liminf_{k\rightarrow\infty}\|Iu_k\|_{L^2(B(0,R))}\leqslant\liminf_{k\rightarrow\infty}\|u_k\|_{L^2(B(0,R))}\leqslant\limsup_{t\rightarrow
T^*}\|u(t)\|_{L^2(B(0,R))},\ee and \be\label{5.47}
m_n\leqslant\liminf_{k\rightarrow\infty}\|In_k\|_{L^1(B(0,R))}\leqslant\liminf_{k\rightarrow\infty}\|n_k\|_{L^1(B(0,R))}\leqslant\limsup_{t\rightarrow
T^*}\|n(t)\|_{L^1(B(0,R))}. \ee

\end{proof}


\end{document}